\RequirePackage{fix-cm}
\PassOptionsToPackage{table}{xcolor}
\PassOptionsToPackage{hyphens}{url}
\documentclass[smallextended]{svjour3}
\usepackage{graphicx}
\usepackage{amsmath,amsfonts,amssymb}

\usepackage{amsthm}
\usepackage{mathtools}
\usepackage{booktabs}
\usepackage[utf8]{inputenc}
\usepackage[T1]{fontenc}
\usepackage{lmodern}
\usepackage{enumitem}
\usepackage{natbib}
\usepackage[colorlinks=true,linkcolor=blue,citecolor=blue,urlcolor=blue]{hyperref}
\usepackage[table]{xcolor}
\usepackage{tikz}
\usetikzlibrary{shapes.geometric,arrows.meta,positioning}

\theoremstyle{plain}

\newtheorem{Prop}{Proposition}
\newtheorem{Lem}{Lemma}
\newtheorem{Cor}{Corollary}
\theoremstyle{definition}

\newtheorem{Rem}{Remark}

\newcommand{\Ical}{\mathcal{I}}
\newcommand{\Tcal}{\mathcal{T}}
\newcommand{\Hcal}{\mathcal{H}}
\newcommand{\Kcal}{\mathcal{K}}
\newcommand{\Rcal}{\mathcal{R}}
\newcommand{\Rv}{\Rcal^{\mathrm v}}
\newcommand{\Rb}{\Rcal^{\mathrm b}}
\newcommand{\Zplus}{\mathbb{Z}_{\ge 0}}
\newcommand{\cg}{c^{\mathrm g}}
\newcommand{\cv}{c^{\mathrm v}}
\newcommand{\cb}{c^{\mathrm b}}
\graphicspath{{figs/}}
\journalname{}
\renewcommand{\makeheadbox}{}
\hypersetup{
  pdftitle={Electric Vehicle Scheduling and Vehicle-to-Grid Integration in Microgrids},
  pdfauthor={Nathan Cho; Andrea Lodi; Anna Scaglione},
  pdfsubject={A route covering and energy scheduling formulation}
}

\begin{document}

\title{Electric Vehicle Scheduling and Vehicle-to-Grid Integration in Microgrids}
\subtitle{A route covering and energy scheduling formulation}
\author{Nathan Cho \and
        Andrea Lodi\thanks{Corresponding author: al748@cornell.edu} \and Anna Scaglione} \authorrunning{N.
        Cho et al.}
\institute{
  Nathan Cho \at
    Operations Research and Information Engineering, Cornell Tech, New York, NY, USA\\
    \email{nc437@cornell.edu} \and Andrea Lodi \at
    Jacobs Technion-Cornell Institute, Cornell Tech and Technion - IIT, New York, NY, USA\\
    \email{al748@cornell.edu} \and Anna Scaglione \at
    Electrical and Computer Engineering, Cornell Tech, New York, NY, USA\\
    \email{as337@cornell.edu} }
\date{}
\maketitle

\begin{abstract}
We jointly schedule electric trucks and stationary batteries in an isolated microgrid with solar and fossil generation.
Truck routes determine when their batteries can charge or return energy through vehicle-to-grid (V2G).
The model assigns every transport task, selects assets, and minimizes daily asset and fossil generation costs.

For a fixed task sequence and charging opportunities, energy scheduling is solved by a continuous network flow model.
We solve the fleet-level linear programming (LP) relaxation by column generation.
A dynamic program jointly generates truck routes and energy profiles using fixed state-of-charge increments.
We give conditions under which this discretization is exact.

The method solves the LP relaxation for instances with up to 1{,}050 tasks.
On the 10--450-task benchmark ladder, every reported integer schedule is within 1.4\% of its LP lower bound.
In matched comparisons, stationary storage provides most of the combined storage and V2G benefit and can determine feasibility under tight generation limits.
Solar surplus relative to task traction energy helps explain investment tradeoffs; charger cost and truck availability affect when bidirectional charging becomes economical.
These results support site-specific investment comparisons.
\end{abstract}

\section{Introduction}\label{sec:intro}

Military operations in remote locations face substantial challenges in securing reliable energy, owing to logistical complexity and high transportation cost.
Supplying remote bases with diesel is costly and exposes operations to vulnerable fuel supply chains \citep{ndupress}.

On-site renewable generation can reduce reliance on delivered fossil fuel \citep{ourenergypolicy}.
When parked, electric vehicles (EVs) with vehicle-to-grid (V2G) capability can also act as mobile storage.
A recent demonstration at a military base combined a bidirectional fleet with on-site solar and stationary storage \citep{lbnl2025miramar}.
Because solar output varies over the day, storage can shift energy from periods of surplus to periods of deficit.
This computational study uses a remote military base as its setting, but the formulation also applies more broadly to sites that coordinate vehicle fleets with isolated microgrids.

We schedule a fleet of electric trucks to serve all transport tasks while minimizing fossil generation and asset costs.
The model extends the Electric Vehicle Scheduling Problem (EVSP) by including state of charge (SoC) dynamics, V2G operation, stationary storage, and a microgrid energy balance with solar and fossil generation.
Existing electric routing and scheduling models account for battery limits and charging decisions \citep{schneider2014electric,keskin2016partial, desaulniers2016exact}.
Recent work also studies V2G for buses under fixed transit schedules or grid-connected demand response \citep{yetkin2024cooptimizing,zhou2025charging}.
Section~\ref{sec:lit} explains how these settings differ from ours.
In our model, routing and charging decisions determine the fossil generation required by the microgrid.
Each truck column represents a feasible sequence of transport tasks together with its charging and discharging schedule.
The number of columns can grow exponentially with the number of tasks, so we solve the LP relaxation by column generation.

\paragraph{Contributions.}
The paper makes the following contributions.
\begin{enumerate}[label=(\roman*),leftmargin=2.4em,itemsep=3pt]
\item We formulate the problem as a route covering model coupled to a microgrid power balance.
The balance determines fossil generation from the selected routes and their energy profiles (Section~\ref{sec:model}).
\item For a fixed route and task sequence, feasibility and optimal energy scheduling are determined by a continuous network flow LP.
We also show that identical stationary batteries can be represented by one aggregate storage block (Sections~\ref{sec:model} and~\ref{sec:structure}).
\item We solve the route pricing problem by dynamic programming over discretized SoC levels.
If charging is 100\% efficient and the battery capacity and rate limit, boundary charge levels, and travel energy requirements are multiples of the chosen step, the discretized and continuous pricing problems have the same optimum.
The final integer solution uses the routes generated during column generation, and we report its gap from the LP lower bound (Section~\ref{sec:cg} and Corollary~\ref{cor:lattice}).
\item We compare four technology configurations to separate the effects of stationary storage and V2G.
We then examine how their value changes with solar surplus, fleet size, task timing, and charger cost (Section~\ref{sec:exp}).
\item A separate feasibility study shows that, under tight generation capacity, some tested task sets require stationary storage when every deployed truck must serve at least one transport task (Section~\ref{ssec:infra}).
\end{enumerate}

This paper is organized as follows.
Sections~\ref{sec:lit}--\ref{sec:notation} review related work, state the assumptions, and introduce the notation.
Sections~\ref{sec:model}--\ref{sec:cg} present the formulation, structural results, and column generation method.
Section~\ref{sec:exp} reports the computational study, and Section~\ref{sec:concl} concludes.


\section{Literature review}\label{sec:lit}

Conventional vehicle scheduling has an established exact and heuristic literature
\citep{bunte2009overview,pepin2009comparison}.
Electric vehicle scheduling now covers several operational features;
\citet{perumal2022electric} provide a review.
Examples include mixed electric and diesel fleets with range and refueling constraints \citep{li2019mixed}, charging during vehicle duties
\citep{van2017scheduling}, partial charging at capacitated stations
\citep{devos2024electric}, robust scheduling under stochastic traffic
\citep{tang2019robust}, and time-of-use tariffs and peak charging loads
\citep{WU2022322}.
\citet{desaulniers2016exact} develop exact branch-price-and-cut algorithms for electric vehicle routing.
Electric vehicle scheduling has also been addressed with a direct integer programming model \citep{janovec2019exact}, joint scheduling and charging formulations \citep{sassi2017electric}, and an exact Benders-type decomposition
\citep{alvo2021exact}.
Battery degradation and nonlinear charging can also affect the resulting schedules \citep{zhang2021optimal}.

A separate body of work studies interactions between electric vehicles and power systems.
\citet{sedighizadeh2020microgrid} coordinate EV charging and discharging with responsive loads in a grid-connected microgrid under uncertainty.
In a marina with on-site solar and stationary storage, bidirectional charging increased solar self-consumption \citep{jozwiak2023optimal}.
A U.S. Army field demonstration connected tactical vehicles with onboard generators and separate battery units in a V2G--V2V microgrid
\citep{masrur2017military}.

The closest scheduling papers differ from ours in either the energy setting or the transportation decisions.
\citet{devos2024electric} solve the LP relaxation of an electric vehicle scheduling model with partial charging and capacitated stations by column generation.
They obtain integer schedules with two heuristics that use column generation and do not allow bidirectional discharge.
\citet{zhou2025charging} jointly schedule charging and V2G discharging for electric buses at multiple depots.
Their buses buy and sell energy on a connected grid at exogenous prices that vary by time, and their method is an adaptive large-neighborhood search.
Their model does not include an isolated microgrid balance or stationary storage optimized jointly with the fleet.

\citet{yetkin2024cooptimizing} jointly optimize an urban power system and an electric transit fleet under renewable uncertainty.
Their transit schedule is fixed, whereas our model assigns required transport tasks to truck duties.
\citet{terada2025multi} size microgrid assets and charging infrastructure and schedule EV charging and V2G under solar, load, and outage scenarios.
Vehicle arrival times, departure times, and energy requirements are inputs;
their model does not assign transport tasks.
Our model instead assigns timed tasks and constructs complete truck duties.
\citet{caustur2025integrated} combine V2G with on-site solar and storage while taking each bus's trip assignment and timetable as fixed.

Related grid-connected studies coordinate electric bus operations with on-site photovoltaic energy and stationary storage
\citep{liu2024electric,bie2024optimal,luke2025optimal}.
Other V2G studies use exogenous tariffs or market prices rather than fossil generation determined by an isolated microgrid balance
\citep{fei2023exploring,zhou2025charging,son2025urban}.

Table~\ref{tab:lit} summarizes these distinctions.
Our contribution is to assign timetabled transport tasks to truck duties while optimizing V2G and microgrid operation.
Fossil generation, on-site solar use, and stationary storage deployment are determined within the same optimization model.

\begin{table}[htbp]
\centering
\caption{Selected comparisons with related studies.
Y $=$ yes, N $=$ no, and P $=$ partial.
Task assignment means that required timetable trips are assigned to vehicle duties.
For \citet{son2025urban}, P denotes aggregate service assignment rather than individual trip chains.
For \citet{luke2025optimal}, P denotes assignment of preconstructed service blocks to buses; the model does not construct the blocks from timetable trips.
PV denotes on-site photovoltaic generation.
Endogenous generation means that dispatchable local generation is determined within the optimization model.
Y under uncertainty denotes a stochastic or robust optimization model.
For this paper, P under uncertainty denotes out-of-sample weather evaluation, not stochastic optimization.}
\label{tab:lit}
\footnotesize
\setlength\tabcolsep{3.5pt}
\begin{tabular}{lcccccc}
\toprule
 & task & & on-site & stat. & endog. & uncer- \\
 & assignment & V2G & PV & storage & gen. & tainty \\
\midrule
\citet{devos2024electric} & Y & N & N & N & N & N \\
\citet{WU2022322} & Y & N & N & N & N & N \\
\citet{zhou2025charging} & Y & Y & N & N & N & N \\
\citet{yetkin2024cooptimizing} & N & Y & N & N & Y & Y \\
\citet{son2025urban} & P & Y & N & N & N & Y \\
\citet{caustur2025integrated} & N & Y & Y & Y & N & N \\
\citet{liu2024electric} & N & N & Y & Y & N & N \\
\citet{bie2024optimal} & Y & N & Y & Y & N & N \\
\citet{luke2025optimal} & P & N & Y & Y & N & N \\
\citet{terada2025multi} & N & Y & Y & Y & Y & Y \\
\midrule
This paper & Y & Y & Y & Y & Y & P \\
\bottomrule
\end{tabular}
\end{table}

\section{Problem setting and assumptions}\label{sec:setting}

A geographically isolated base must serve a set of transport tasks with a fleet of electric trucks while minimizing fossil generation and asset costs.
The microgrid has a known base electricity demand and on-site solar generation; trucks draw power to charge and can return power to the grid (V2G) or to one another (V2V).
We make the following assumptions.

\begin{enumerate}[label=(A\arabic*),leftmargin=2.4em,itemsep=2pt]
\item \emph{Forecast.}
Demand and solar generation are known over the horizon; the model is deterministic.
Section~\ref{ssec:limits} discusses stochastic extensions.
\item \emph{Battery limit for each vehicle.}
Each truck has a finite battery capacity $G$, obeys SoC limits, and has a charge/discharge rate $\rho$.

\item \emph{Charging capacity.}
In each time block, total charging from the microgrid to trucks and stationary storage is bounded by an exogenous system wide limit.
The associated scarcity price enters the column generation pricing problem.

\item \emph{Generation capacity.}
Fossil generation in each block is bounded by the generator capacity, $g_t\le\bar g_t$.

\item \emph{Linear charging and energy loss.}
Charging and discharging are linear.
If the grid supplies $q$ units of energy during charging, the battery stores $(1-\eta)q$, where $\eta\in[0,1)$ is the charging loss fraction.
When $\eta=0$, charging is 100\% efficient.
The model also includes a small throughput cost $\sigma$, defined in Section~\ref{sec:model}, to favor schedules with less battery cycling.
The model does not explicitly represent battery degradation or nonlinear charging curves \citep{montoya2017electric,zhang2021optimal}.
V2G cycling can affect battery life \citep{peterson2010lithium,uddin2017possibility}.
Appendix~\ref{app:gammadetail} tests an additional cost on discharge from both trucks and stationary batteries.
Its upper value, \$0.13/kWh, was chosen to be close in magnitude to a 2030 mean compensation estimate of EUR~0.132/kWh that includes both degradation and V2G infrastructure \citep{sagaria2025vehicle}; it is intended as a stress test, not an estimate of cycling wear alone.
\citet{preis2023battery} report much smaller degradation costs in their base cases.

\item \emph{Fleet sizing.}
The numbers of trucks and stationary batteries are decision variables.
Each deployed asset has a positive cost, so an optimal solution uses finitely many assets.
\end{enumerate}

\section{Notation}\label{sec:notation}

Table~\ref{tab:notation} summarizes the notation.

\begin{table}[htbp]
\centering
\caption{Notation.}\label{tab:notation}
\footnotesize
\setlength{\tabcolsep}{5pt}
\renewcommand{\arraystretch}{1.1}
\begin{tabular}{lp{0.58\linewidth}}
\toprule
\multicolumn{2}{l}{\textbf{Sets and indices}}\\
$i,j\in\Ical$ & transport tasks (trips); $n=|\Ical|$ \\
$t\in\Tcal,\quad \Tcal=\{1,\dots,T\}$ & time blocks over the planning horizon \\
$h\in\Hcal,\quad \Hcal^0=\Hcal\cup\{o\}$ & charging station $h$; depot $o$; all charging locations \\
$k\in\Kcal,\quad \Kcal=\Kcal^{\mathrm v}\cup\Kcal^{\mathrm b}$ & assets: trucks and stationary batteries \\
$r\in\Rcal,\quad \Rcal=\Rv\cup\Rb$ & feasible columns: truck routes and battery schedules; $\Rcal'\subseteq\Rcal$ is the restricted set \\
\midrule
\multicolumn{2}{l}{\textbf{Data}}\\
$[\tau^{\mathrm s}_i,\tau^{\mathrm e}_i]$, $\ell^{\mathrm s}_i,\ell^{\mathrm e}_i$ & start and end times (fixed) and start/end locations of task $i$ \\
$\varepsilon_i$ & traction energy required to perform task $i$ \\
$d_{ij}$ & deadhead distance from task $i$ to task $j$; a deadhead is empty travel between tasks, with proportional time and energy \\
$G,\ \rho,\ \eta$ & battery capacity; charge/discharge rate per block; charging loss fraction \\
$D_t,\ P_t$ & base electricity demand and solar (PV) generation in block $t$ \\
$\Delta_t=D_t-P_t$ & net demand; $\Delta_t^+=\max\{0,\Delta_t\}$, $\Delta_t^-=\max\{0,-\Delta_t\}$ \\
$\cg$ & cost per unit of fossil generation \\
$\bar g_t$ & fossil generation capacity in block $t$ ($=\infty$ if uncapped) \\
$\sigma$ & small throughput cost per unit of gross charging or discharging \\
$\cv,\ \cb$ & fixed deployment cost of a truck / a battery \\
$a_{ir}\in\{0,1\}$ & $1$ iff column $r$ covers task $i$ \\
$e_{rht}\in\mathbb R$ & net grid energy of column $r$ at location $h\in\Hcal^0$ in block $t$ \\
$e^+_{rht}:=\max\{e_{rht},0\}$ & positive part of the net grid draw by column $r$ at location $h$ in block $t$ \\
$\kappa_t$ & system wide charging capacity in block $t$ \\
$L$ & number of physical (task) locations; each may host a charging station \\
\midrule
\multicolumn{2}{l}{\textbf{Variables and dual prices}}\\
$x_r\in\Zplus$ & number of times column $r$ is used \\
$g_t\ge 0$ & fossil generation dispatched in block $t$ \\
$\alpha_i,\ \mu_t,\ \nu_t$ & coverage dual ($\alpha_i$ free), power balance dual $\mu_t\ge0$, and nonnegative charging capacity scarcity price $\nu_t$ \\
\bottomrule
\end{tabular}
\end{table}

The coefficient $e_{rht}$ is signed: $e_{rht}>0$ means that column $r$ draws energy from the grid at location $h$ in block $t$, while $e_{rht}<0$ means that it returns energy to the grid.
Lemma~\ref{lem:nosim} shows that an optimal schedule can be chosen that never charges and discharges simultaneously.
For such a schedule,
$e^+_{rht}=\max\{e_{rht},0\}$
is the energy drawn for charging and is the quantity counted against the system wide charging capacity.
Thus, the power balance uses the net interaction $e_{rht}$, while the charging capacity constraint uses $e^+_{rht}$.

\section{The unified EVSP--V2G formulation}\label{sec:model}

\paragraph{Route and battery columns.}
A \emph{column} $r$ is a feasible plan for one asset over the horizon.
A truck column $r\in\Rv$ specifies its tasks, deadhead movements, visits to charging stations, and feasible energy profile.
The coefficient $a_{ir}$ records whether the column covers task $i$.
A stationary battery column $r\in\Rb$ covers no transport task and specifies only an energy profile.
Truck and battery columns interact with the microgrid through the same signed coefficients $e_{rht}$.
Their deployment cost components are $\cv$ and $\cb$, respectively; the complete column cost is defined below.
Because both types enter the same power balance through $e_{rht}$, the master chooses mobile and stationary storage jointly.

\paragraph{Master problem.}
Let $x_r$ count the use of column $r$ and let $g_t\ge 0$ be the fossil generation in block $t$.
The microgrid EVSP--V2G is
\begin{subequations}\label{eq:master}
\begin{align}
\min_{x,g}\quad & \sum_{t\in\Tcal}\cg\,g_t \;+\; \sum_{r\in\Rcal} c_r\,x_r \label{eq:obj}\\
\text{s.t.}\quad
& \sum_{r\in\Rv} a_{ir}\,x_r = 1, && \forall i\in\Ical, \label{eq:cover}\\
& g_t \;-\; \sum_{r\in\Rcal}\sum_{h\in\Hcal^0} e_{rht}\,x_r \;\ge\; \Delta_t, && \forall t\in\Tcal, \label{eq:balance}\\
& \sum_{r\in\Rcal}\sum_{h\in\Hcal^0} e^+_{rht}\,x_r
  \;\le\; \kappa_t, && \forall t\in\Tcal, \label{eq:chgcap}\\
& 0 \le g_t \;\le\; \bar g_t, && \forall t\in\Tcal, \label{eq:gpos}\\
& x_r \in \Zplus, && \forall r\in\Rcal. \label{eq:int}
\end{align}
\end{subequations}
The objective~\eqref{eq:obj} includes fossil generation and asset costs.
Constraint~\eqref{eq:cover} serves every task exactly once.
Constraint~\eqref{eq:balance} requires fossil generation and net storage discharge to meet net demand $\Delta_t$.
Its slack represents unused energy.
We assume that unused solar can be curtailed and that any other excess can be safely absorbed, for example by a dump load.
A site without this capability would need a finite limit on excess absorption.
Constraint~\eqref{eq:chgcap} limits total charging in each time block.
Constraint~\eqref{eq:gpos} imposes nonnegativity and the generation limit; $\bar g_t=\infty$ represents an uncapped block.
Constraint~\eqref{eq:int} requires an integer number of columns.
For a truck route covering at least one task, the equality constraints~\eqref{eq:cover} already force $x_r\le1$.
Using the route twice would cover each of its tasks twice.
No explicit upper bound is therefore needed for truck routes.
We retain the general integer form because stationary battery schedules may be used with multiplicity.

\paragraph{Fossil generation and the positive part relation.}
For a fixed feasible selection of columns, $\cg>0$ implies that an optimal solution uses the smallest value of fossil generation allowed by the power balance constraint.
Hence,
\[
g_t=
\max\left\{
0,\,
\Delta_t+\sum_{r,h}e_{rht}x_r
\right\}.
\]
This value must not exceed $\bar g_t$; otherwise the selected columns are infeasible.
The equality follows from the linear power balance constraint and the positive generation cost.
Charging increases fossil generation only after the expression inside the maximum becomes positive.
Holding all other grid interactions fixed, discharging reduces fossil generation while $g_t>0$ and has no generation value after $g_t$ reaches zero.

The formulation allows net storage discharge $\sum_{r,h}(-e_{rht})x_r$ to exceed the contemporaneous deficit $\Delta_t^+$.
Discharge beyond this level cannot reduce fossil generation further, but allowing it prevents the SoC discretization from excluding otherwise useful columns.
For example, if the deficit is $30$\,kWh and the SoC step is $50$\,kWh, the only useful discretized profile may discharge $50$\,kWh.
That profile discharges 20\,kWh more than the deficit, but the excess does not reduce the objective.

\begin{Rem}[Economic dispatch interpretation]\label{rem:dispatch}
Keeping $g_t$ as an explicit variable in the objective is a modeling choice.
The dual $\mu_t$ of the power balance is the marginal value of energy in block $t$.
Without a binding generation or fuel constraint, $0\le\mu_t\le\cg$.
At a block with positive generation and no binding generation limit, $\mu_t=\cg$.
A binding generation or fuel constraint can add a scarcity rent $\lambda_t\ge0$, so that $\mu_t=\cg+\lambda_t$ at a block with positive generation.
Sections~\ref{sec:cg} and~\ref{ssec:resil} return to this case.
These prices encourage charging when prices are low and discharging when prices are high.

Eliminating $g_t$ would replace its linear cost by
\[
\cg\sum_t\left[\Delta_t+\sum_{r,h}e_{rht}x_r\right]^+,
\]
together with the corresponding generation cap condition.
The resulting formulation has the same optimum, but it hides the marginal price.
Keeping $g_t$ also makes it easier to add standard power system constraints, including multiple generators, convex generation costs, ramp limits, and operating reserves.
We therefore keep generation cost $\cg\sum_tg_t$ separate from asset costs.
\end{Rem}

\begin{Rem}[Charging infrastructure]\label{rem:infra}
Constraint~\eqref{eq:chgcap} limits total charging energy in each time block.
It applies to both trucks and stationary batteries.
We write $\nu_t\ge0$ for the marginal cost of charging capacity in block $t$.

The base model does not limit the number of plugs.
Such a limit could be added through a constraint on column occupancy.
The model limits each asset's discharge rate but does not impose a common limit on the total energy returned by all batteries in one block.
Such a constraint can be added if the microgrid has a binding limit on returned power.

The experiments do not model electrical network constraints among charging locations.
Separate limits for each location can be added by replacing \eqref{eq:chgcap} with one constraint for every location and time block.
In the reported experiments, finite charging caps are used only for depot charging; the experiments with multiple stations leave the charging cap unbounded.
Thus, the experiments do not test separate charging caps at each station.
\end{Rem}

\begin{Rem}[Vehicle-to-vehicle transfer and energy cycling]\label{rem:v2v}
All charging locations in the experiments are connected to one electrical node.
At that node, one truck discharging while another truck charges is already represented by their two signed grid interactions.
A separate V2V variable would not change the power balance.
Such a variable would become relevant in a network with line limits or with a V2V efficiency different from a transfer through the grid.

The common node can produce alternative optimal solutions that circulate energy without changing fossil generation.
When $\eta=0$, a cycle between trucks loses no energy.
Even when $\eta>0$, a cycle can have no cost if it occurs in blocks with $\mu_t=\nu_t=0$.
We use two parameters to address these cases.
The physical parameter $\eta$ accounts for charging loss.
For each column $r$, let $c_{rht}\ge0$ and $d_{rht}\ge0$ denote the energy drawn for charging and returned through discharging at location $h$ in block $t$, so that
\[
e_{rht}=c_{rht}-d_{rht}.
\]
Define
\[
C_r=\sum_{h\in\Hcal^0}\sum_{t\in\Tcal}c_{rht},
\qquad
D_r=\sum_{h\in\Hcal^0}\sum_{t\in\Tcal}d_{rht}.
\]
The throughput term in the column formulation is
\[
\sigma\sum_{r\in\Rcal}(C_r+D_r)x_r.
\]
This term favors schedules with less cycling.
Lemma~\ref{lem:nosim} shows that, because $\sigma>0$, a minimum-cost energy schedule cannot charge and discharge simultaneously.

We set $\sigma=0.025$ per 100\,kWh model unit, $1/1{,}600$ of the \$0.40/kWh generation cost.
The term discourages unnecessary cycling, and because its coefficient is finite, it can also change the ordering of solutions whose other costs are very close.
The sensitivity test in Appendix~\ref{app:gammadetail} separately applies an additional cost to discharge from both trucks and stationary batteries.
\end{Rem}

\paragraph{Costs.}
The generation term $\cg\sum_t g_t$ contains the cost of energy.
Using the definitions above, the column cost is
\[
c_r=c^{\mathrm{asset}}_r+\sigma(C_r+D_r),
\qquad c^{\mathrm{asset}}_r\in\{\cv,\cb\}.
\]
When charging and discharging do not occur simultaneously, gross throughput equals $\sum_{h,t}|e_{rht}|$.
The implementation uses this form for truck columns.
The aggregate battery block in Proposition~\ref{rem:agg} applies the same $\sigma$ directly to its gross charge and discharge variables.
Charging loss affects the energy dynamics, while the resulting generation cost is accounted for through $g_t$.

\begin{Prop}[Exact aggregation of identical stationary storage]\label{rem:agg}
For real $k\ge0$, let $\mathcal P_k$ be the set of feasible charge, discharge, and state schedules in the linear battery model with capacity $kG$ and rate $k\rho$, allowing simultaneous charging and discharging.
Assume that the batteries are connected to the same electrical node and have linear dynamics and costs.
Use a boundary condition that scales with $k$: either $s_T=s_0$, or $s_0=kG$ with an unrestricted terminal state.
Then $\mathcal P_0=\{0\}$ and $\mathcal P_k=k\mathcal P_1$ for every $k>0$, with costs scaling by $k$.
For an integer $k>0$,
\[
\mathcal P_k=\left\{\sum_{j=1}^{k}p^j:p^j\in\mathcal P_1\right\}.
\]
Thus, $k$ identical batteries are equivalent to one aggregate block with capacity $kG$ and rate $k\rho$.
The statement for real $k$ describes the LP relaxation, in which the battery count is continuous.
\end{Prop}
\begin{proof}
For $k>0$, all constraints scale linearly with capacity and rate, so $\mathcal P_k=k\mathcal P_1$.
For every positive integer $k$, the sum of $k$ feasible unit schedules satisfies the aggregate dynamics and bounds.
Conversely, any $p\in\mathcal P_k$ can be disaggregated by assigning the schedule $p/k\in\mathcal P_1$ to each unit.
\end{proof}

\paragraph{Implementation of stationary storage.}
Proposition~\ref{rem:agg} shows when several stationary batteries can be represented exactly by a single storage block.
The result requires identical batteries connected to the same modeled electrical node, with linear dynamics and costs and the same boundary condition.
These assumptions hold for the stationary batteries in our computational model, so we use one aggregate block.
The variable $N_{\mathrm b}$ counts the installed batteries; it is continuous in the LP relaxation and integer in the final master problem.
In each time block, aggregate charging and discharging are bounded by $N_{\mathrm b}\rho$, and aggregate SoC is bounded by $N_{\mathrm b}G$.
The corresponding storage dynamics and the applicable boundary condition are included directly in every restricted master problem.
Column generation is used only for truck routes.
This representation may fail when batteries have different initial or terminal requirements, unit-specific connection or operating constraints, or nonlinear or unit-specific degradation costs.
We retain the battery column notation $\Rb$ in the general formulation because it can represent such heterogeneous batteries.

\section{Structural results}\label{sec:structure}

This section presents a network flow formulation for fixed route energy scheduling and a submodularity result for fixed energy profiles.
Proofs not given in the text appear in Appendix~\ref{app:proofs}.

\subsection{Energy scheduling for a fixed route}

A fixed route skeleton specifies the task sequence, deadhead movements, and the blocks and locations in which the truck is parked and connected.
It does not assume that the route satisfies the energy constraints.
The LP below tests energy feasibility and, when the route is feasible, finds its energy profile of minimum cost.

For a fixed skeleton, build a time-expanded SoC network.
Storage nodes record the battery state over time, and carryover arcs have capacity $G$.
A charging arc draws energy from the microgrid and stores the fraction $1-\eta$.
A discharge arc returns stored energy to the microgrid.
Fixed withdrawals represent task and deadhead energy.
Charging and discharging arcs are available only while the asset is parked and connected, and their capacities impose the rate limit $\rho$.

\begin{Prop}[Energy scheduling as a network flow LP]\label{prop:flow}
For any fixed candidate skeleton, the time-expanded energy flow LP represents the linear SoC dynamics, rate limits, and boundary conditions, allowing simultaneous charging and discharging.
The candidate route satisfies the energy constraints if and only if this LP is feasible.
For the nonnegative energy and charging-capacity prices and nonnegative throughput cost used in pricing, an optimal LP schedule can be chosen without simultaneous charging and discharging and attains the minimum cost over schedules with that restriction.
When $\eta=0$, this is an ordinary network flow LP.
Its network incidence matrix is totally unimodular, so, when the LP is feasible, integral data admit an integral optimal flow.
When $\eta>0$, it is a generalized flow LP.
In both cases, charging and discharging are continuous decisions.
\end{Prop}

The piecewise-linear charging formulations in \citet{montoya2017electric} and \citet{froger2019improved} use binary variables to select charging segments.
Under the linear charging assumption used here, the network flow formulation above remains continuous.

The charging capacity constraint uses $e^+_{rht}=\max\{e_{rht},0\}$.
This quantity equals gross charging when the asset does not charge and discharge simultaneously.
The following lemma establishes that an optimal schedule can be chosen with this property.

\begin{Lem}[No simultaneous charge and discharge]\label{lem:nosim}
Let $c_t$ and $d_t$ denote the grid energy drawn for charging and returned by discharging in block $t$, so $e_t=c_t-d_t$.
Fix a skeleton, prices $\mu_t\ge0$ and $\nu_t\ge0$, and a throughput cost $\sigma\ge0$.
Assume that the fixed-skeleton energy scheduling problem is feasible.
At least one energy schedule of minimum cost does not charge and discharge simultaneously.
When $\sigma>0$, every energy schedule of minimum cost has this property.
For any such schedule, $e_t^+=c_t$.
\end{Lem}
\begin{proof}
Let
\[
B=\{t\in\Tcal:c_t>0,\ d_t>0\}.
\]
For each $t\in B$, set
\[
\epsilon_t=\min\{c_t,d_t/(1-\eta)\}
\]
and replace $(c_t,d_t)$ by
\[
\bigl(c'_t,d'_t\bigr)
=\bigl(c_t-\epsilon_t,\ d_t-(1-\eta)\epsilon_t\bigr).
\]
For $t\notin B$, leave $c_t$ and $d_t$ unchanged.
For every $t\in B$,
\[
(1-\eta)c'_t-d'_t
=(1-\eta)c_t-d_t,
\]
so the complete SoC path is unchanged.
The modified quantities are nonnegative, and either $c'_t=0$ or $d'_t=0$.
Moreover,
\[
e'_t=c'_t-d'_t=e_t-\eta\epsilon_t\le e_t.
\]
Thus neither the energy price term nor the charging capacity term can increase, while gross throughput decreases by $(2-\eta)\epsilon_t$ in every block in $B$.
Because $\mu_t,\nu_t,\sigma\ge0$, the modified schedule costs no more than the original schedule.
If $B\ne\varnothing$ and $\sigma>0$, it costs strictly less.
Hence, when $\sigma=0$, some minimum-cost schedule has no simultaneous charging and discharging; when $\sigma>0$, no minimum-cost schedule can contain such a block.
\end{proof}

\begin{Rem}[Boundary conditions in the pricing problem]
The boundary condition on SoC is a modeling choice.
In Proposition~\ref{prop:flow}, it changes only the source and sink conditions.
In the pricing dynamic program, it changes the initial labels and the terminal acceptance rule.
Section~\ref{sec:cg} states the conventions used, and Section~\ref{ssec:repro} measures what is at stake between them.
\end{Rem}

\paragraph{Residual energy redispatch.}
Appendix~\ref{app:proofs} gives a separate polymatroid result for a simplified redispatch problem in which the routes and a feasible energy schedule are fixed.

\subsection{The fixed profile displacement value is submodular}

We next consider a simpler selection problem after the truck routes have been fixed.
Each candidate profile $p$ can deliver a fixed amount $b_{p,t}$ of surplus energy to deficit block $t$.
The following conditions separate this problem from the full route covering model.

\paragraph{Fixed profile conditions.}
The truck routes are fixed, and the candidate profiles form a finite set $\mathcal F$.
For each $p\in\mathcal F$ and $t\in\Tcal$, the value $b_{p,t}\ge0$ is fixed.
Selecting a profile does not change the feasibility or deliverability of another profile.
The profiles interact only through the useful displacement limit $\Delta_t^+$ in each block.

\begin{Lem}[Submodularity of fixed profile displacement]\label{thm:submod}
Under the fixed profile conditions above, define
\[
V(S)=\sum_{t\in\Tcal}\min\Big\{\sum_{p\in S} b_{p,t},\ \Delta_t^+\Big\},
\qquad S\subseteq\mathcal{F}.
\]
Then, $V$ is monotone and submodular on $2^{\mathcal{F}}$.
\end{Lem}
\begin{proof}
Fix $t\in\Tcal$ and write $m_t(S)=\sum_{p\in S}b_{p,t}$.
Since $b_{p,t}\ge0$, this is a nonnegative modular set function on $2^{\mathcal{F}}$.
The block transform $\Phi_t(z)=\min\{z,\Delta_t^+\}$ is the pointwise minimum of the affine maps $z\mapsto z$ and $z\mapsto\Delta_t^+$, hence concave, and it is nondecreasing.
For $A\subseteq B$ and $p\notin B$, we have $m_t(A)\le m_t(B)$.
Concavity of $\Phi_t$ therefore gives
\[
\Phi_t(m_t(A)+b_{p,t})-\Phi_t(m_t(A))
\ge
\Phi_t(m_t(B)+b_{p,t})-\Phi_t(m_t(B)).
\]
Thus, the contribution from block $t$ is monotone and submodular.
The sum over all blocks preserves both properties, so $V$ is monotone and submodular.
\end{proof}

\begin{Rem}[Scope and use of Lemma~\ref{thm:submod}]\label{rem:submodscope}
Lemma~\ref{thm:submod} applies after the routes and energy delivery profiles have been fixed.
Taking expectations preserves the result when every scenario uses the same feasible profile set and satisfies the same conditions.
The result need not hold when route covering, a physical V2V rendezvous, or shared constraints make the profiles interact.
\end{Rem}

\subsection{Limitations of the formulation}\label{ssec:limits}

The model treats net demand $\Delta_t=D_t-P_t$ as known.
In practice, both demand and solar generation are uncertain.
The deterministic formulation can be solved again as forecasts change, but a stochastic model would represent this uncertainty directly.

Section~\ref{ssec:commit} takes an empirical first step toward uncertainty by fixing schedules for a planning solar profile and evaluating their dispatch over one year of realized weather data.

\section{Solution method}\label{sec:cg}
\subsection{Route set and computational complexity}

The number of feasible route columns is exponential in the number of tasks and time blocks, so an explicit master over all columns does not scale.
The classical single-depot vehicle scheduling problem can be solved as an assignment problem \citep{bertossi1987matching,freling2001models}.
Adding battery constraints changes the complexity.
The single-depot electric vehicle scheduling problem is NP-hard even without V2G \citep{li2014transit,sassi2017electric,perumal2022electric}.
Appendix~\ref{app:proofs} gives a self-contained reduction showing strong NP-hardness for a restriction of our formulation.

\begin{Prop}[Strong NP-hardness with energy constraints]\label{thm:np}
The single-depot restriction of \eqref{eq:master} with energy constraints, and therefore the full EVSP--V2G formulation, is strongly NP-hard.
This holds with one vehicle type, no V2G, no stationary storage, $\eta=0$, and temporally disjoint tasks at the same location.
\end{Prop}

This reduction locates the hardness in battery capacity rather than spatial routing.
In the reduction, all tasks are at the same location and charging is unavailable until the final recharge period.
Assigning tasks to trucks is therefore a bin-packing decision: the tasks' energy requirements must fit within each truck's battery capacity.
Appendix~\ref{app:proofs} gives the full reduction from 3-Partition.

The pricing subproblem contains 0--1 knapsack: it selects task rewards subject to an energy limit.
The dynamic program with discretized SoC described below is consequently pseudo-polynomial.
Strong NP-hardness rules out a pseudo-polynomial exact algorithm for the full problem unless P$=$NP.

\subsection{Column generation method}
We solve the LP relaxation of the master problem by column generation.
Starting from a small set of routes, each iteration solves the restricted master LP and uses its dual prices to find a route with negative reduced cost.
Such a route is added to the restricted master, and the process repeats.
Figure~\ref{fig:algorithm-flowchart} summarizes this procedure.
To initialize the restricted master, we temporarily add one artificial column for each task.
It covers only that task, uses no energy, and carries a large penalty.
These columns make the initial LP feasible whenever the microgrid without the fleet is feasible.
\begin{figure}[htbp]
\centering
\resizebox{\linewidth}{!}{%
\begin{tikzpicture}[
node distance=1.45cm,
every node/.style={font=\small},
block/.style={
rectangle, draw, rounded corners,
align=center,
minimum width=5.0cm,
minimum height=0.9cm
},
decision/.style={
diamond, draw,
align=center,
aspect=2.4,
inner sep=1pt,
minimum width=4.2cm,
minimum height=1.25cm
},
line/.style={->, thick}
]

\node[block] (start)
{Initialize the restricted master with available truck\\
routes and temporary artificial columns};

\node[block, below of=start] (rmp)
{Solve restricted master LP\\
and obtain dual prices $(\alpha,\mu,\nu)$};

\node[block, below of=rmp, minimum height=1.3cm] (pricing)
{Use dynamic programming to find the truck route\\
with the lowest reduced cost};

\node[decision, below of=pricing, yshift=-0.7cm] (test)
{Negative reduced cost\\
route found?};

\node[block, below of=test, xshift=-2.9cm, yshift=-1.15cm] (add)
{Add generated columns\\
to restricted master};

\node[block, below of=test, xshift=2.9cm, yshift=-1.15cm] (complete)
{LP column generation\\
complete};

\draw[line] (start) -- (rmp);
\draw[line] (rmp) -- (pricing);
\draw[line] (pricing) -- (test);

\draw[line] (test) -- node[above left] {yes} (add);

\draw[line]
(add.west) -- ++(-1.0,0)
|- (rmp.west);

\draw[line] (test) -- node[above right] {no} (complete);

\end{tikzpicture}}

\caption{Overview of the column generation procedure used to solve the master LP relaxation.}
\label{fig:algorithm-flowchart}
\end{figure}

\paragraph{Restricted master and duals.}
Solve the LP relaxation of \eqref{eq:master} (replace \eqref{eq:int} by $x_r\ge0$) over a restricted set $\Rcal'$.
Let $\alpha_i$ be the free duals of coverage~\eqref{eq:cover}, and let $\mu_t\ge0$ be the duals of power balance~\eqref{eq:balance}.
The conventional dual of the charging capacity constraint~\eqref{eq:chgcap} is nonpositive; we write it as $-\nu_t$, where $\nu_t\ge0$ is the corresponding scarcity price.
The dual $\mu_t\ge0$ is the marginal value of relaxing the power balance requirement in block $t$.
Without a binding generation capacity constraint, pricing of $g_t$ gives $\mu_t\le\cg$; when the cap $g_t\le\bar g_t$ is present and binding, or when an aggregate fuel budget $\sum_t g_t\le B$ binds (Section \ref{ssec:resil}), $\mu_t$ can also include the scarcity value of limited generation.

\paragraph{Reduced costs and pricing.}
The reduced cost of a candidate column $r$ is
\begin{equation}\label{eq:rc}
\bar c_r \;=\; c_r \;-\; \sum_{i\in\Ical} a_{ir}\,\alpha_i
\;+\; \sum_{h\in\Hcal^0}\sum_{t\in\Tcal}\mu_t\, e_{rht}
\;+\; \sum_{h\in\Hcal^0}\sum_{t\in\Tcal}\nu_t\, e^+_{rht}.
\end{equation}
The term involving $\mu_t$ charges or credits the route for its net exchange with the microgrid.
It therefore favors charging when marginal energy is inexpensive and discharging when generation is costly or scarce.
The term involving $\nu_t$ penalizes charging in a block whose charging capacity is scarce.

For a fixed sequence of tasks and station visits, the remaining charging and discharging problem is the continuous network flow LP of Proposition~\ref{prop:flow}.
Conceptually, pricing therefore contains two parts:
\begin{itemize}[leftmargin=1.4em,itemsep=2pt]
\item a routing decision that chooses the tasks, their sequence, and the station visits; and
\item an energy scheduling decision that minimizes the energy terms in \eqref{eq:rc}, together with the small activity cost $\sigma$, subject to the SoC dynamics.
\end{itemize}
The implementation combines these decisions in one labeling dynamic program, rather than solving a separate LP at every label.
Related labeling methods appear in exact electric vehicle routing algorithms with linear charging \citep{desaulniers2016exact}, in exact algorithms for charging on a fixed route with nonlinear charging \citep{froger2019improved}, and in charge scheduling with time-dependent energy prices \citep{klein2023electric}.
The direct arc costs used here are given below.
Stationary storage does not require a pricing subproblem because it is represented directly by the aggregate master block of Proposition~\ref{rem:agg}.

\paragraph{Initial route set.}
The first restricted master contains one feasible route serving only task $i$ for each task for which such a route can be constructed.
It also contains one temporary artificial column for every task.
For a route serving one task, the truck recharges at the depot after completing the task whenever the remaining time is sufficient to restore its terminal SoC.
These recharge amounts are continuous and need not be multiples of $\delta$.

Before solving the first restricted master LP, we call the pricing dynamic program repeatedly to construct additional routes that serve several tasks.
We call this enriched initial route set the warm start.
In each call, tasks not yet included in a previous warm start route receive a positive provisional reward, while tasks already included receive a large negative reward.
The returned route is added, and the procedure continues with the tasks not yet included.
The artificial columns provide temporary coverage while column generation searches for real routes.

The aggregate stationary storage block is present from the first restricted master and needs no warm start.
Across 32 matched runs, the change from a cold start to a warm start averaged 3.4\% fewer iterations, ranging from 10.9\% fewer to 2.1\% more.
We use the warm start mainly to obtain useful initial dual prices.

\medskip

\paragraph{Pricing dynamic program.}
Truck routing and energy scheduling are solved together by a label-setting dynamic program on an acyclic graph indexed by time, location, and SoC.
A label has state $(t,\ell,s,k)$, where $t$ is the time block, $\ell$ is a physical location, $s$ is the truck's SoC, and $k$ records whether the route has covered at least one task.

The dynamic program represents SoC at the levels
\(
   \{0,\delta,2\delta,\ldots,G\},
\)
where $\delta$ is the SoC step.
All steps used in the reported experiments divide $G$, so the number of levels is
\(
n_s=G/\delta+1.
\)
The principal computational studies use $\delta=25$\,kWh, which gives 29 levels for a 700\,kWh truck battery.
The coarser studies that use $\delta=50$\,kWh are identified in Section~\ref{sec:exp}.

A generated truck route enters the master through its task incidence vector and fixed energy profile.
The master then optimizes continuous fossil generation and aggregate stationary storage dispatch over the available truck routes.

The flag $k$ must equal one at a terminal label, so every generated truck route serves at least one task.
When stationary storage can be purchased with the same capacity, rate, location, and charging loss as a truck, no tighter boundary condition, and cost $\cb\le\cv$, a truck used only for storage is weakly dominated.
These conditions hold in the reference configuration, but need not hold in every cost or battery-capacity sensitivity test.
When stationary storage is unavailable or these dominance conditions do not hold, requiring every truck to serve a task defines the fleet considered in the experiment.
Allowing additional parked trucks could only improve those objective values, so the reported values are conservative with respect to that extension.

\paragraph{Boundary conditions.}
The reported experiments use two boundary conventions for truck SoC:
\begin{itemize}[leftmargin=1.4em,itemsep=2pt]
\item \emph{Full recharge}: $s_0=s_T=G$.
Both the initial and terminal SoC are fixed at full capacity.
This is the default readiness requirement for trucks.
\item \emph{Periodic}: $s_0=s_T$.
The initial and terminal SoC must be equal, but the model chooses their common value.
\end{itemize}
Periodic operation therefore includes full recharge as one option but may choose a lower repeated level.
For example, starting below capacity leaves room to absorb early solar surplus without increasing fossil generation.
The aggregate stationary battery is periodic under both truck conventions.
Because the truck boundary condition changes which routes are feasible, columns are generated separately for each convention.
A separate free start diagnostic, in which trucks begin full at no cost and may finish at any SoC, is reported only in Section~\ref{ssec:repro}.

The graph contains wait, charge, discharge, deadhead, and task arcs.
Charge and discharge arcs change $s$ by one or more steps of $\delta$.
Deadhead and task arcs consume traction energy, and task arcs also contribute the coverage reward $-\alpha_i$.
Every arc advances time.
The graph is therefore acyclic, and a single forward sweep obtains shortest paths even when some reduced cost arc weights are negative.

\paragraph{When the SoC discretization is exact.}
The integrality statement in Proposition~\ref{prop:flow} gives the following condition under which pricing with discretized SoC agrees with the continuous LP for a fixed route.

\begin{Cor}[Exactness of the SoC discretization when $\eta=0$]
\label{cor:lattice}
Suppose that charging incurs no energy loss, so $\eta=0$.
Assume that $G$, $\rho$, every task and deadhead energy, and every fixed boundary SoC are integer multiples of $\delta$.
For full recharge or periodic operation, the pricing dynamic program returns the optimum of the continuous fixed route energy model.
If the dynamic program and the continuous model consider the same route skeletons, their master LPs also have the same value.
\end{Cor}

\begin{proof}
Fix a route and the master dual prices.
By Lemma~\ref{lem:nosim}, an optimal energy schedule can be chosen without simultaneous charging and discharging, so $e_t^+=c_t$.
Replacing the term $\nu_t e_t^+$ by the linear term $\nu_t c_t$ does not change the optimum: the two terms agree for that schedule, and $c_t\ge e_t^+$ for every schedule.

Now divide every energy quantity by $\delta$.
Proposition~\ref{prop:flow} gives an ordinary minimum cost flow problem with integral capacities and balances.
Under full recharge, the boundary level is fixed at $G/\delta$.
Under periodic operation, an arc from the final storage node to the initial storage node represents the common boundary level and has capacity $G/\delta$.
Flow integrality therefore gives an optimum at the SoC levels used by the dynamic program.

Thus the continuous and discretized pricing problems have the same value for each route skeleton.
The flow feasible set for each fixed route is a bounded polytope, so it has finitely many extreme points.
Since pricing agrees for every master dual vector, LP duality gives equality of the two master LP values.
\end{proof}

The main experiments use $\delta=25$\,kWh.
This step divides the 50\,kWh task and deadhead energies and the 175\,kWh limit per block, so the corollary applies when $\eta=0$.
In the studies with a 175\,kWh block limit and $\delta=50$\,kWh, the pricing dynamic program can use at most 150\,kWh per block.
The structured timing benchmarks instead use a 50\,kWh block limit and a matching 50\,kWh step, so the corollary applies to those lossless instances as well.
Appendix~\ref{app:extra} reports the effect of refining this approximation.

The result does not apply when $\eta>0$, because the corresponding generalized flow matrix need not be integral.
When $\eta=0$, the charging capacity term remains linear: by Lemma~\ref{lem:nosim}, with $\sigma>0$ an optimum does not charge and discharge simultaneously, so $e^+_{rht}$ is the gross charging flow.

The discretization check provides a numerical comparison.
Across 16 matched instances with $\eta=0$, the LP values obtained with $\delta=25$ and $\delta=12.5$\,kWh agree to every recorded digit.
When $\eta=0.15$, the corollary no longer applies.
Across the four matched cases, refining the SoC step from 12.5 to 6.25\,kWh lowers the LP objective by between 0.1\% and 2.0\%.
\begin{Prop}[Pricing and finite convergence for the discretized column family]\label{prop:cg}
(i) Every arc of the pricing graph advances time.
A forward sweep therefore finds a minimum reduced cost column among all columns admitted by the SoC discretization, even when some arc costs are negative.

(ii) Let $\mathcal C$ be the finite union of the initial real columns and the columns searched by the pricing dynamic program.
Suppose that artificial columns make the initial restricted master feasible.
Each continuing iteration adds a previously absent column with negative reduced cost under the current restricted master duals.
The algorithm stops only after exact pricing under those duals finds no omitted negative reduced cost column in $\mathcal C$.
It then terminates finitely at an optimum of the master with artificial columns over $\mathcal C$.
If the solution uses no artificial column, it is also optimal for the master without artificial columns over $\mathcal C$.
Otherwise, a separate feasibility test is required.
\end{Prop}
\begin{proof}
(i) Process nodes in increasing time.
Whenever a node is processed, all its predecessors already have their shortest-path labels.
This argument uses acyclicity and does not require nonnegative arc costs.

(ii) Since $\mathcal C$ is finite and every continuing iteration adds a previously absent column, only finitely many additions are possible.
At termination, restricted master optimality gives the dual inequalities for generated columns, and exact pricing gives them for every omitted column in $\mathcal C$.
The current dual is therefore feasible for the full master with artificial columns.
The restricted master primal solution has the same objective value, so strong duality proves optimality.
If that solution uses no artificial column, it is feasible and optimal for the master without artificial columns.
Otherwise, the augmented optimum is only a lower bound for that master.
\end{proof}

\paragraph{Implemented column family.}
The column family $\mathcal C$ includes the finite warm start routes and the routes searched by the discretized dynamic program.
Some warm start routes use recharge amounts between the dynamic program's SoC levels.
Termination is checked by pricing at the current master duals; stabilized duals are used only to propose additional columns.

\paragraph{Phase I feasibility test.}
Artificial columns are temporary coverage placeholders rather than physical truck routes.
Because their variables are continuous in the LP and their penalty is finite, a small positive amount can remain even when a solution using only real columns exists.
Residual artificial coverage in the economic LP therefore triggers a separate feasibility test.

If artificial coverage remains after pricing at the current restricted master dual vector finds no route with reduced cost below $-\epsilon$, we solve an auxiliary LP over the same column family $\mathcal C$.
Its objective is the total artificial coverage; the original route, storage, and generation costs are ignored \citep{desrosiers2026branch}.
Column generation continues until pricing at the current restricted master dual vector finds no column with reduced cost below $-\epsilon$.
If the optimum exceeds the tolerance bound described below, the LP without artificial columns is infeasible, so no integer solution over $\mathcal C$ exists.
If the optimum is zero, that LP is feasible, although integer feasibility has not yet been established.
In that case, we retain the columns found during Phase I and resume column generation with the original cost objective.
We accept a solution only if it uses no artificial coverage; otherwise the outcome remains unresolved.
For instances with $\eta=0$ that satisfy Corollary~\ref{cor:lattice}, the same feasibility conclusion holds for the continuous profile master over the same family of route skeletons.

The implementation stops pricing when the smallest reduced cost is at least $-\epsilon$.
Decreasing every coverage dual $\alpha_i$ by $\epsilon$ raises the reduced cost of every real route by at least $\epsilon$, because every such route serves at least one task.
The adjusted dual is feasible for the full master, and its objective is at most $n\epsilon$ below the restricted master value.
We therefore use $z_{\mathrm{RMP}}-n\epsilon$ as the lower bound for the full master.
A positive Phase I objective certifies infeasibility only when it exceeds $n\epsilon$.

An artificial column that appears only in an integer solve that reaches the time limit has a different meaning: the solver did not find an incumbent over the generated routes that uses no artificial columns.
Such an outcome is reported as unresolved rather than infeasible.

Routes are elementary because each task is available only at its unique start time and every arc advances time; a path cannot cover the same task twice.

\paragraph{Arc costs for energy actions.}
The dynamic program assigns the energy terms in~\eqref{eq:rc} directly to its arcs.
Suppose a charging arc at location $h$ and time $t$ draws $q$ kWh from the microgrid.
Then $e_{rht}=e^+_{rht}=q$, and the arc contributes $(\mu_t+\nu_t)q$ to reduced cost.
Because only the fraction $1-\eta$ enters the battery, increasing SoC by $\delta$ requires $q=\delta/(1-\eta)$.
A discharge arc that supplies $q$ kWh to the microgrid has $e_{rht}=-q$ and $e^+_{rht}=0$, and therefore contributes $-\mu_tq$.
A charge or discharge arc carrying $q$ kWh also contributes $\sigma q$.

\paragraph{Complexity.}
Recall that $n=\lvert\Ical\rvert$ is the number of tasks, $T=\lvert\Tcal\rvert$ is the number of time blocks, $L$ is the number of physical locations, $H=\lvert\Hcal^0\rvert$ is the number of charging locations, and $n_s$ is the number of SoC levels.
Let
\[
u^{+}=\lfloor(1-\eta)\rho/\delta\rfloor,
\qquad
u^{-}=\lfloor\rho/\delta\rfloor
\]
be the maximum numbers of charging and discharging steps in one block.
For a fixed initial SoC, the graph has $O(TLn_s)$ states; the binary flag $k$ changes this count only by a constant factor.
Deadhead transitions contribute $O(TL^2n_s)$ arcs, energy transitions contribute $O(THn_s(u^{+}+u^{-}))$, and task transitions contribute $O(nn_s)$.
The implementation scans the task list across the $T$ blocks, which adds $O(Tn)$ work.
Since the graph is acyclic, one forward sweep therefore costs
\[
O\!\bigl(
Tn_s\{L^2+H(u^{+}+u^{-})\}+nn_s+Tn
\bigr).
\]
Full recharge and a periodic boundary fixed at a specified level each use one such sweep.
When the common periodic level is chosen by the model, the implementation runs one sweep for each of the $n_s$ possible levels.
Its worst-case bound is therefore $n_s$ times the expression above.

Since $n_s=\Theta(G/\delta)$ and $u^\pm=\Theta(\rho/\delta)$, the sweep over fixed levels is pseudo-polynomial in the scaled quantities $G/\delta$ and $\rho/\delta$.
Its energy transition term is quadratic in $1/\delta$; choosing the periodic level adds another factor of $n_s$.
The implementation evaluates the SoC levels within each sweep by vectorized operations.

For the reported instances, one pricing call takes a few hundredths of a second.
Column generation requires a few tens of iterations on the planning instances, a few hundred at intermediate sizes, and 1{,}400--1{,}700 on the 1{,}000-task instances in Table~\ref{tab:diagbig}.

\paragraph{Dominance.}
At each state $(t,\ell,s,k)$, representing the time block, location, SoC, and whether the route has served a task, the dynamic program keeps only the label with minimum cost.
Labels at different SoC levels are retained separately.
A label with high SoC is not assumed to dominate one with low SoC, because the latter may have more room to absorb inexpensive energy later in the day.
Explicit enumeration of the 29 SoC levels is inexpensive here.
Resource-constrained labeling often uses stronger dominance rules that compare both cost and resource consumption \citep{desaulniers2016exact}; such rules may be useful with a finer SoC step or additional resources.

\paragraph{Dual smoothing.}
To reduce large changes in the dual prices between iterations, the first pricing call uses the current restricted master dual vector.
From the second iteration onward, the first pricing call uses the average of the current and previous restricted master dual vectors,
\(
  \hat d=\tfrac12 d^{\mathrm{prev}}+\tfrac12d,
  \qquad
  d=(\alpha,\mu,\nu).
\)
Every route found using $\hat d$ is evaluated with the current duals and is added only if its current reduced cost is below the stopping threshold.
If this call adds no route, pricing is repeated with the current duals.
We declare LP column generation converged only when pricing at the current duals finds no route with reduced cost below that threshold.
After the column generation loop, additional pricing calls with reproducibly perturbed duals broaden the route set used in the final integer solve.

\paragraph{Final integer solution.}
After LP column generation and the additional pricing calls, we fix the generated columns and solve the resulting restricted master as a mixed integer program.
An incumbent that uses no artificial columns is feasible for \eqref{eq:master}.
Since only the generated columns are considered, this heuristic solve does not certify integer optimality over the entire column family $\mathcal C$.
Section~\ref{ssec:method} explains the reported gap.

\section{Computational study}\label{sec:exp}

\subsection{Experimental setup}\label{ssec:setup}

All experiments use a 24-hour horizon divided into 48 half-hour periods.
Unless a sensitivity analysis changes these values, truck batteries have capacity $G=700$\,kWh and a charge and discharge rate of 350\,kW.
This rate corresponds to $\rho=175$\,kWh per period.
The principal experiments use an SoC step of $\delta=25$\,kWh.
The structured timing benchmarks in Table~\ref{tab:diag} and their extensions in Fig.~\ref{fig:scal} use a 100\,kW charging rate, a 50\,kWh SoC step, and a generation cost of \$0.05/kWh.
The randomly generated task instances in Table~\ref{tab:diagbig} also use a 50\,kWh step.
Appendix~\ref{app:extra} reports the discretization checks and the design of the planning experiments.

Restricted master LPs were solved through SciPy's HiGHS interface.
CBC solved the final MIPs over generated routes for the structured benchmark family, and Gurobi 12.0.3 solved them in the remaining computational studies.
The 10--450-task benchmark timing runs in Fig.~\ref{fig:scal} were performed on an Apple M5 Pro system with 48\,GB of memory.
The extended benchmark runs used single nodes of a Linux cluster with Intel Xeon E5-2600~v3 processors (2.3--2.4\,GHz).

The benchmark family uses demand and solar profile shapes adapted from \citet{almeida2020energy} and scaled for a synthetic coastal-island microgrid.
For context, a San Nicolas Island report gives a peak load of approximately 975\,kW \citep{maynard2019sni}; the benchmark is not calibrated to reproduce that load.
Before rounding net demand to the model's energy units, the reference profiles contain approximately 27.8\,MWh of daily demand and 14.7\,MWh of daily solar generation.
After rounding, the daily net deficit and surplus are 16.9 and 3.7\,MWh, respectively.
A $k\times$ solar case multiplies the reference solar output in every time block by $k$, while demand is unchanged.

Two-hour tasks require 200\,kWh.
The benchmark family is used for the timing and boundary condition experiments.
The planning family uses the same demand shapes and the economic parameters in Table~\ref{tab:prov}; the table is placed in Appendix~\ref{app:extra}.
A second scalability family uses one-hour tasks with randomly sampled origins, destinations, and start times on maps with 4 or 15 task locations.
The model still chooses the number of trucks; only the task sets are randomly generated.

Four technology configurations are compared: charge-only trucks without stationary storage, charge-only trucks with stationary storage, V2G trucks without stationary storage, and V2G trucks with stationary storage.
In tables and figures, battery energy storage system (BESS) denotes the stationary storage option.

The 270-instance technology decomposition and the main experiments on generation capacity, charging capacity, annual weather, and bidirectional charger premiums use the same generated routes for matched configurations with the same truck technology.
The V2G route set also contains every charge-only route.
Other sensitivity studies construct route sets separately for each scenario.
Appendix~\ref{app:extra} states which construction is used.

\subsection{Computational performance}\label{ssec:method}

The computations have two stages.
Column generation first solves the LP relaxation to the stated pricing tolerance.
The generated routes are then fixed, and the resulting restricted master is solved as a MIP.
This second stage produces a feasible integer incumbent, but no additional routes are generated.
Optimality in this stage therefore applies only to the fixed generated route set.

\begin{figure}[htbp!]
\centering
\includegraphics[width=\linewidth]{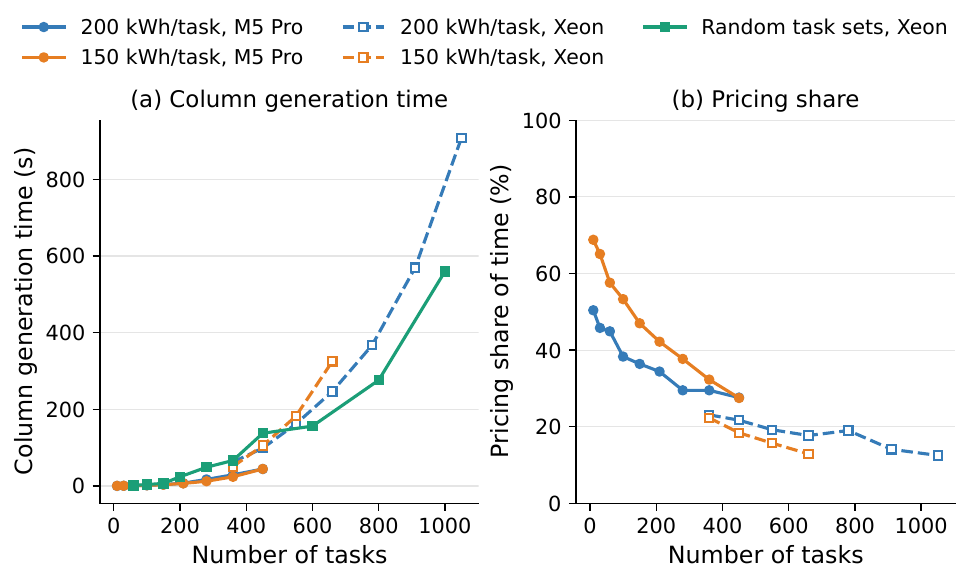}
\caption{Column generation time (a) and the share of time spent in pricing (b).
Solid curves show the 10--450-task benchmark timings on the Apple M5 Pro system.
Dashed curves with open markers show the same instances solved again and extended on Intel Xeon cluster nodes.
At the overlapping sizes, the reported LP objectives agree to the reported precision.
The two timing series are shown separately because the processors differ.
Filled squares show the randomly generated task instances on the cluster.
Table~\ref{tab:diagbig} reports the randomly generated task instances with 4 and 15 task locations.}
\label{fig:scal}
\end{figure}

The structured benchmark runs in Table~\ref{tab:diag} are short, and the reported time includes column generation and the final MIP over generated routes.
A 100-task run takes less than two seconds.
At 450 tasks, the 200\,kWh and 150\,kWh settings take 44.6 and 49.1 seconds and have LP-to-incumbent gaps of 0.03\% and 1.17\%, respectively.
For the randomly generated instances with charging only at the depot in Fig.~\ref{fig:scal}, column generation reaches the pricing tolerance for 1{,}000 tasks in less than fifteen minutes.
When charging is available at the 4 or 15 task locations, Table~\ref{tab:diagbig} reports median column generation times of 3.6--4.6 hours for 1{,}000 tasks.
The final MIP then has a separate time limit of 1{,}800 seconds.
For both 100-task rows, this MIP is optimal over the generated routes; every larger row reaches the time limit.
At 1{,}000 tasks, the LP-to-incumbent gaps are 41--44\%.

Additional scalability tests use the same benchmark family with task energy requirements of 150 and 200\,kWh.
For the 200\,kWh instances, reported gaps remain below 0.12\% through 1{,}050 tasks.
Among the completed 150\,kWh runs, the gaps reach about 2.5\% beyond 500 tasks.
The 150\,kWh cases above 660 tasks did not finish in the available time and are not reported.

Pricing accounts for about half of the column generation time on the smallest instances, less than one third at 450 tasks, and at most 2\% at 1{,}000 tasks.
Thus the restricted master LP, rather than the pricing DP, is the main computational bottleneck at large sizes.
The final integer restricted master is also difficult on the randomly generated task instances.

Let $z_{\mathrm{RMP}}$ be the final restricted master LP value.
The pricing tolerance gives the lower bound
\[
\mathrm{LB}=z_{\mathrm{RMP}}-n\epsilon.
\]
This correction is below the displayed precision, so the tables report $z_{\mathrm{RMP}}$ as the LP bound.
The bound applies to the implemented column family $\mathcal C$.
It also applies to the continuous-profile master over the same route skeletons when Corollary~\ref{cor:lattice} applies; otherwise it need not be a lower bound for that larger model.
For an integer incumbent with value $z^{\mathrm I}$, the reported gap is
\[
100\frac{z^{\mathrm I}-z_{\mathrm{RMP}}}{z^{\mathrm I}}.
\]
This measure is not the optimality gap reported by the MIP solver.

The LP-to-incumbent gap has up to three sources: the LP integrality gap, the use of only the generated routes in the final MIP, and unfinished MIP search at the time limit.
When the final MIP is marked ``r.opt.,'' only the first two sources remain.
We do not estimate these components separately.
The reported gap therefore does not certify integer optimality over the full route family.

\begin{table}[htbp]
\centering
\caption{Solver diagnostics on the 10--450-task ladder, whose maps grow with the instance ($L=2$--$10$ task locations).
The LP bound and MIP incumbent are in dollars per day; columns $=$ truck columns generated; pricing \% $=$ share of column generation time spent in the pricing DP; total time includes column generation and the final MIP.
Negative gaps with absolute value below 0.01\% are due to solver rounding and are reported as 0.}
\label{tab:diag}
\footnotesize
\setlength\tabcolsep{3.8pt}
\begin{tabular}{rrrrrrrrr}
\toprule
$n$ & $L$ & LP bound & MIP incumbent & gap \% & columns & CG iters & pricing \% & total (s) \\
\midrule
\multicolumn{9}{l}{\emph{each task energy 200\,kWh}}\\
10 & 2 & 1{,}060.60 & 1{,}060.60 & 0.00 & 60 & 22 & 50 & 0.2 \\
30 & 3 & 1{,}571.80 & 1{,}571.80 & 0.00 & 103 & 37 & 46 & 0.4 \\
60 & 4 & 2{,}463.60 & 2{,}463.60 & 0.00 & 150 & 42 & 45 & 0.6 \\
100 & 5 & 3{,}806.25 & 3{,}806.25 & 0.00 & 266 & 102 & 38 & 1.7 \\
150 & 6 & 5{,}519.62 & 5{,}569.65 & 0.90 & 399 & 185 & 36 & 3.7 \\
210 & 7 & 7{,}603.72 & 7{,}608.75 & 0.07 & 561 & 287 & 34 & 7.4 \\
280 & 8 & 10{,}058.55 & 10{,}058.55 & 0.00 & 836 & 492 & 30 & 17.2 \\
360 & 9 & 12{,}884.10 & 12{,}889.42 & 0.04 & 1{,}124 & 700 & 30 & 29.5 \\
450 & 10 & 16{,}095.45 & 16{,}100.48 & 0.03 & 1{,}354 & 840 & 28 & 44.6 \\
\midrule
\multicolumn{9}{l}{\emph{each task energy 150\,kWh}}\\
10 & 2 & 1{,}035.48 & 1{,}035.47 & 0.00 & 138 & 25 & 69 & 0.4 \\
30 & 3 & 1{,}496.42 & 1{,}496.42 & 0.00 & 175 & 34 & 65 & 0.6 \\
60 & 4 & 2{,}312.85 & 2{,}312.85 & 0.00 & 239 & 56 & 58 & 1.0 \\
100 & 5 & 3{,}454.95 & 3{,}454.95 & 0.00 & 319 & 80 & 53 & 1.6 \\
150 & 6 & 4{,}892.63 & 4{,}892.62 & 0.00 & 449 & 160 & 47 & 3.4 \\
210 & 7 & 6{,}625.88 & 6{,}626.12 & 0.00 & 614 & 265 & 42 & 6.8 \\
280 & 8 & 8{,}654.70 & 8{,}760.82 & 1.21 & 810 & 391 & 38 & 13.1 \\
360 & 9 & 10{,}979.10 & 11{,}134.82 & 1.40 & 1{,}072 & 573 & 32 & 25.8 \\
450 & 10 & 13{,}599.07 & 13{,}759.90 & 1.17 & 1{,}378 & 789 & 28 & 49.1 \\
\bottomrule
\end{tabular}
\end{table}

\begin{table}[htbp]
\centering
\caption{LP scalability and integer incumbents for randomly generated task instances with 4 or 15 task locations; entries are medians over two task draws.
Each column generation run terminates at the stated pricing tolerance using a 50\,kWh SoC step.
The final MIP over generated routes has a time limit of 1{,}800 seconds. ``r.opt.'' means that this MIP was solved to optimality; it does not establish integer optimality over the full route family. ``t.l.'' denotes a feasible incumbent available when the time limit was reached.
The reported LP-to-incumbent gap uses the definition in the text above.
The gap is the median of the gaps for the two task draws; it is not recomputed from the displayed median objectives.
The pricing percentage uses the median pricing and column generation times.}
\label{tab:diagbig}
\scriptsize
\setlength\tabcolsep{2.5pt}
\begin{tabular}{rrrrrrrrrl}
\toprule
$n$ & $L$ & LP bound & MIP incumbent & gap \% & CG iters & pricing \% & CG (s) & MIP (s) & status \\
\midrule
100 & 4 & 4{,}973.8 & 5{,}570.4 & 10.7 & 300 & 20 & 19 & 32 & r.opt. \\
200 & 4 & 10{,}155.9 & 12{,}269.5 & 16.1 & 648 & 6 & 191 & 1{,}800 & t.l. \\
400 & 4 & 19{,}771.7 & 25{,}799.4 & 20.5 & 813 & 2 & 722 & 1{,}800 & t.l. \\
600 & 4 & 29{,}030.9 & 47{,}853.2 & 39.3 & 1{,}002 & $<$1 & 1{,}512 & 1{,}801 & t.l. \\
800 & 4 & 38{,}878.3 & 65{,}076.7 & 40.3 & 1{,}367 & $<$1 & 5{,}580 & 1{,}801 & t.l. \\
1{,}000 & 4 & 48{,}198.4 & 82{,}326.4 & 41.4 & 1{,}720 & $<$1 & 16{,}437 & 1{,}802 & t.l. \\
\midrule
100 & 15 & 6{,}318.4 & 6{,}573.2 & 3.9 & 456 & 46 & 68 & 153 & r.opt. \\
200 & 15 & 12{,}032.0 & 14{,}450.3 & 16.5 & 629 & 25 & 178 & 1{,}801 & t.l. \\
400 & 15 & 23{,}344.6 & 34{,}443.8 & 32.2 & 1{,}004 & 8 & 725 & 1{,}800 & t.l. \\
600 & 15 & 33{,}471.9 & 53{,}473.0 & 37.4 & 1{,}282 & 3 & 2{,}571 & 1{,}801 & t.l. \\
800 & 15 & 44{,}177.0 & 74{,}787.4 & 40.9 & 1{,}468 & 2 & 6{,}676 & 1{,}801 & t.l. \\
1{,}000 & 15 & 54{,}611.6 & 96{,}601.6 & 43.5 & 1{,}724 & $<$1 & 13{,}042 & 1{,}802 & t.l. \\
\bottomrule
\end{tabular}
\end{table}

The 600--1{,}000-task rows test the scalability of the LP column generation method.
Their integer incumbents are not used in the planning conclusions below.

\subsection{Boundary conventions and fossil accounting}\label{ssec:repro}

The reported experiments use full recharge or periodic operation.
Under full recharge, each truck satisfies $s_0=s_T=G$.
Under periodic operation, $s_0=s_T$, but the model chooses the common level.
The aggregate stationary battery is periodic under both conventions.
Periodic operation may begin below full capacity, leaving room for early solar charging while still replacing all energy used over the horizon.
Unless stated otherwise, the truck results use full recharge.

The boundary condition determines which energy schedules can accompany a truck route, so columns are generated separately for each convention.
Table~\ref{tab:boundaryperiodic} in Appendix~\ref{app:boundarydetail} compares periodic operation with an optimized repeated level against five fixed boundary levels.
Relative to the optimized repeated level, full recharge raises median incumbent cost by 14.7\% for charge-only trucks without BESS and 28.1\% for V2G trucks without BESS.
The increases are 2.8\% and approximately 0.1\%, respectively, when BESS is available.
A fixed boundary at half capacity comes within a median 0.5--3.2\% of periodic operation with an optimized repeated level in the three configurations other than V2G with BESS; for V2G with BESS, the median difference is less than 0.1\%.
Thus full recharge imposes a measurable readiness cost.
In these instances, stationary storage makes that restriction much less costly.

The free start convention can also report net export that is financed by unreplaced initial energy.
This export disappears under full recharge on the benchmark instances.
Under either full recharge or periodic operation, summing the storage balances over the day shows that the system formed by the trucks and BESS is a net consumer because it must replace truck traction energy and, when $\eta>0$, charging losses.

\begin{figure}[htbp!]
\centering
\includegraphics[width=\linewidth]{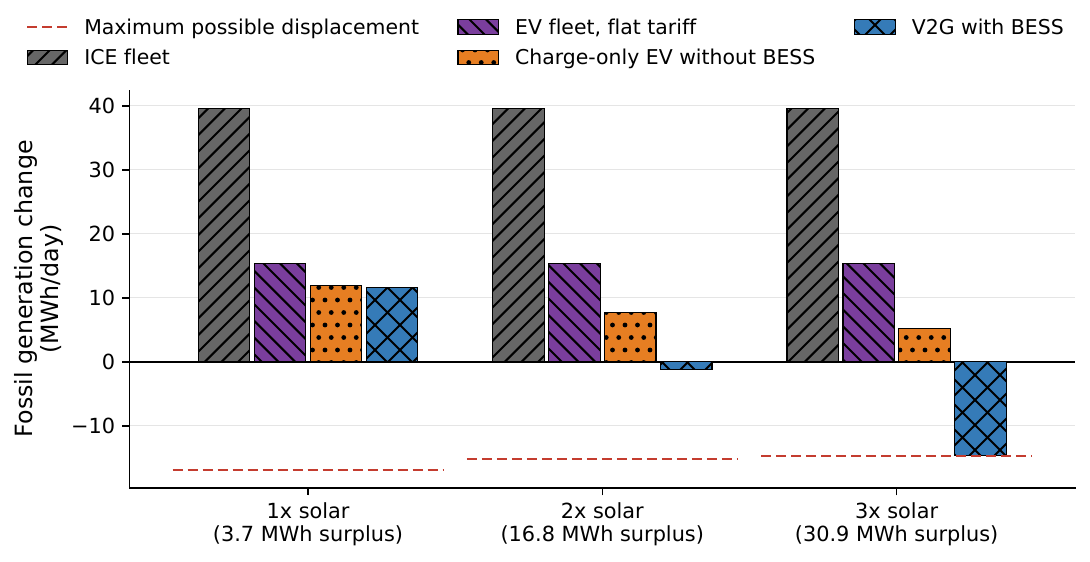}
\caption{Change in fossil generation relative to the microgrid without trucks under full recharge.
The internal combustion engine (ICE) comparison uses 3.3 times the traction energy of an electric truck.
Dashed red segments show the most fossil generation that can be displaced at each solar level.
The flat tariff electric vehicle schedule separates the efficiency benefit of electrification from the benefit of coordinating charging with the microgrid.}
\label{fig:regime}
\end{figure}

Figure~\ref{fig:regime} therefore reports the change in fossil generation relative to the same microgrid operated without trucks.
With little solar, electric trucks add fossil demand.
As solar capacity increases, charging first removes that additional demand and V2G can then displace part of the base's fossil generation.
Table~\ref{tab:export} in Appendix~\ref{app:boundarydetail} gives the corresponding values by workload and solar surplus.
For comparison, the EV benchmark with a flat tariff charges every kWh drawn by the trucks at the constant price of fossil generation and does not coordinate charging with the microgrid.
It isolates the efficiency benefit of electrification.

\subsection{Combined value and its decomposition}\label{ssec:when}

\begin{figure}[htbp!]
\centering
\includegraphics[width=0.9\linewidth]{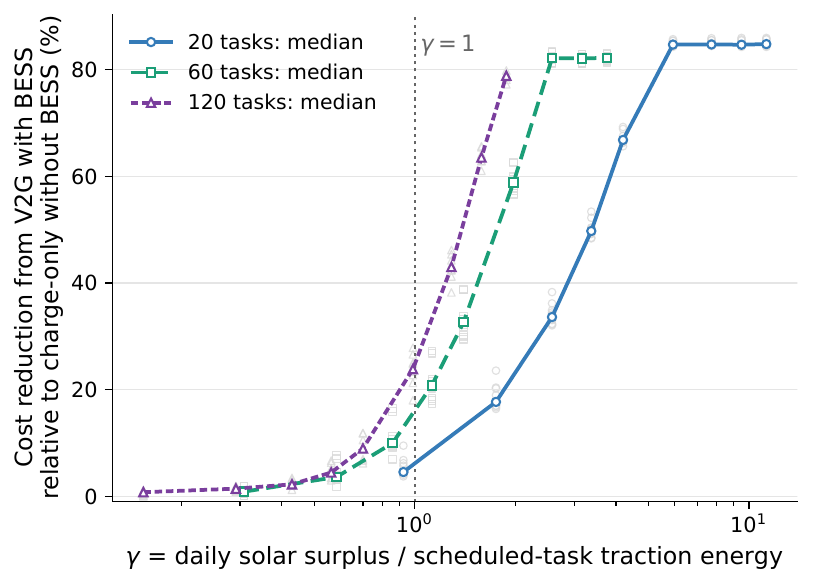}
\caption{Percentage reduction in the final restricted master incumbent cost from allowing V2G and purchasable BESS, relative to charge-only trucks without BESS, in 270 matched instances.
Points show individual task sets; lines show the median for each workload.
Within an instance, the two charge-only configurations share one generated route set; the two V2G configurations share a larger set that also contains the charge-only routes.
No premium for bidirectional chargers is included.
The dotted line marks $\gamma=1$.}
\label{fig:collapse}
\end{figure}

Figure~\ref{fig:collapse} compares charge-only trucks without BESS with V2G trucks that may purchase BESS.
It uses the 270 matched instances in the experiment with four configurations described below: three workloads, nine solar levels, and ten task sets for each workload and solar level.
Each point represents the same base and task set under the two technology configurations.
The horizontal coordinate is
\[
\gamma
=\frac{\sum_{t\in\Tcal}\Delta_t^-}
       {\sum_{i\in\Ical}\varepsilon_i}
=\frac{\text{daily solar surplus}}
       {\text{traction energy for scheduled tasks}}.
\]
The denominator excludes deadhead energy, which depends on the routes the model selects.
The vertical coordinate is the percentage reduction in total daily cost when the model moves from charge-only trucks without BESS to V2G trucks with purchasable BESS.
Positive values mean that the latter configuration is cheaper.

At low values of $\gamma$, daily solar surplus is small relative to task traction energy.
As $\gamma$ increases, more surplus may be available for the batteries to shift to periods that would otherwise require fossil generation.
The separate workload curves in Figure~\ref{fig:collapse} show that $\gamma$ organizes the pattern but does not determine the percentage saving by itself.

Within the experiments that scale the fleet and microgrid together, the percentage cost reduction is largely organized by $\gamma$, but the crossing changes with the cost and operating assumptions.
At a bidirectional charger premium of \$8 per truck per day, the median saving changes sign between $\gamma=0.40$ and $0.45$.
At a \$4 premium, the crossing is near $\gamma=0.35$.

Changing task duration alters when trucks are available to charge and discharge.
Holding task end times fixed and varying duration from one to four hours moves the median crossing between $\gamma=0.32$ and $0.43$, although not monotonically.
Thus, both the amount and timing of the surplus matter.
The ratio $\gamma$ does not include deadhead energy or show whether trucks are connected when the surplus occurs.
The same value of $\gamma$ also gives different percentage savings when base demand is held fixed and the fleet changes size, because the fleet then accounts for a different share of total cost.
The ratio therefore organizes these experiments without defining a universal threshold.
Appendix~\ref{app:gammadetail} gives the individual task sets and sensitivity analyses behind this conclusion.

The four technology configurations identify which resource produces the cost reduction.
For each of three task counts and nine multipliers of solar capacity, we generated ten task sets, giving 270 matched base instances.
All four configurations were solved on each instance.
Based on the final restricted master incumbents, stationary storage accounts for a median 96\% of the combined cost reduction.
The incremental effect of allowing V2G when BESS is available is smaller at low $\gamma$ and larger at high $\gamma$.
The incumbent for V2G with BESS costs less than the incumbent for charge-only trucks with BESS in 265 of the 270 cases; the other five tie.
The LP and incumbent bounds certify a strict reduction in 236 cases, while 34 remain inconclusive.
Appendix~\ref{app:designdetail} gives the interval used for this test.
The component shares need not sum to 100\% because the two resources substitute for one another.

Table~\ref{tab:decomp} summarizes these separate comparisons; its rows are not additive because V2G and stationary storage can substitute for one another.
\begin{table}[htbp]
\centering
\caption{Median cost changes among the four technology configurations over 270 matched instances.
When a share of the combined reduction is reported, its denominator is the reduction from charge-only trucks without BESS to V2G trucks with BESS.}
\label{tab:decomp}
\small
\begin{tabular}{p{0.42\linewidth}p{0.16\linewidth}p{0.28\linewidth}}
\toprule
comparison & $\gamma$ range & median result \\
\midrule
BESS added to charge-only trucks & all & 96\% of combined reduction \\
V2G enabled without BESS & all & 19.6\% of daily cost; 62\% of combined reduction \\
V2G with BESS versus charge-only with BESS & $\gamma\leq0.35$ & 0.6\% of daily cost \\
V2G with BESS versus charge-only with BESS & $\gamma>2$ & 9.3\% of daily cost \\
\bottomrule
\end{tabular}
\end{table}

Appendix Figure~\ref{fig:boundary} reports a separate experiment that varies workload and solar capacity.
At a fixed solar capacity, additional fossil displacement eventually falls as the workload grows and task traction uses more of the available surplus.

\subsection{Infrastructure limits at work}\label{ssec:infra}

The generation cap and the charging cap restrict different parts of the energy balance.
In the experiments, the generation cap multiplier sets $\bar g_t$ in Constraint~\eqref{eq:gpos} equal to that multiplier times
\[
\max_{\tau\in\Tcal}\Delta_\tau^+,
\]
the largest positive net demand of the base before trucks or stationary batteries are added.
Thus, a multiplier of $1.0$ equals the base's peak fossil generation requirement.
The charging cap multiplier similarly sets the common cap $\kappa_t$ in Constraint~\eqref{eq:chgcap} equal to that multiplier times $\max_{\tau\in\Tcal}\Delta_\tau^-$, the largest solar surplus in any block.
This cap limits the combined charging of trucks and stationary batteries in each block.
Table~\ref{tab:design} lists the multipliers used, and Appendix~\ref{app:capdetail} reports the corresponding feasibility and utilization results.

\begin{figure}[htbp!]
\centering
\includegraphics[width=\linewidth]{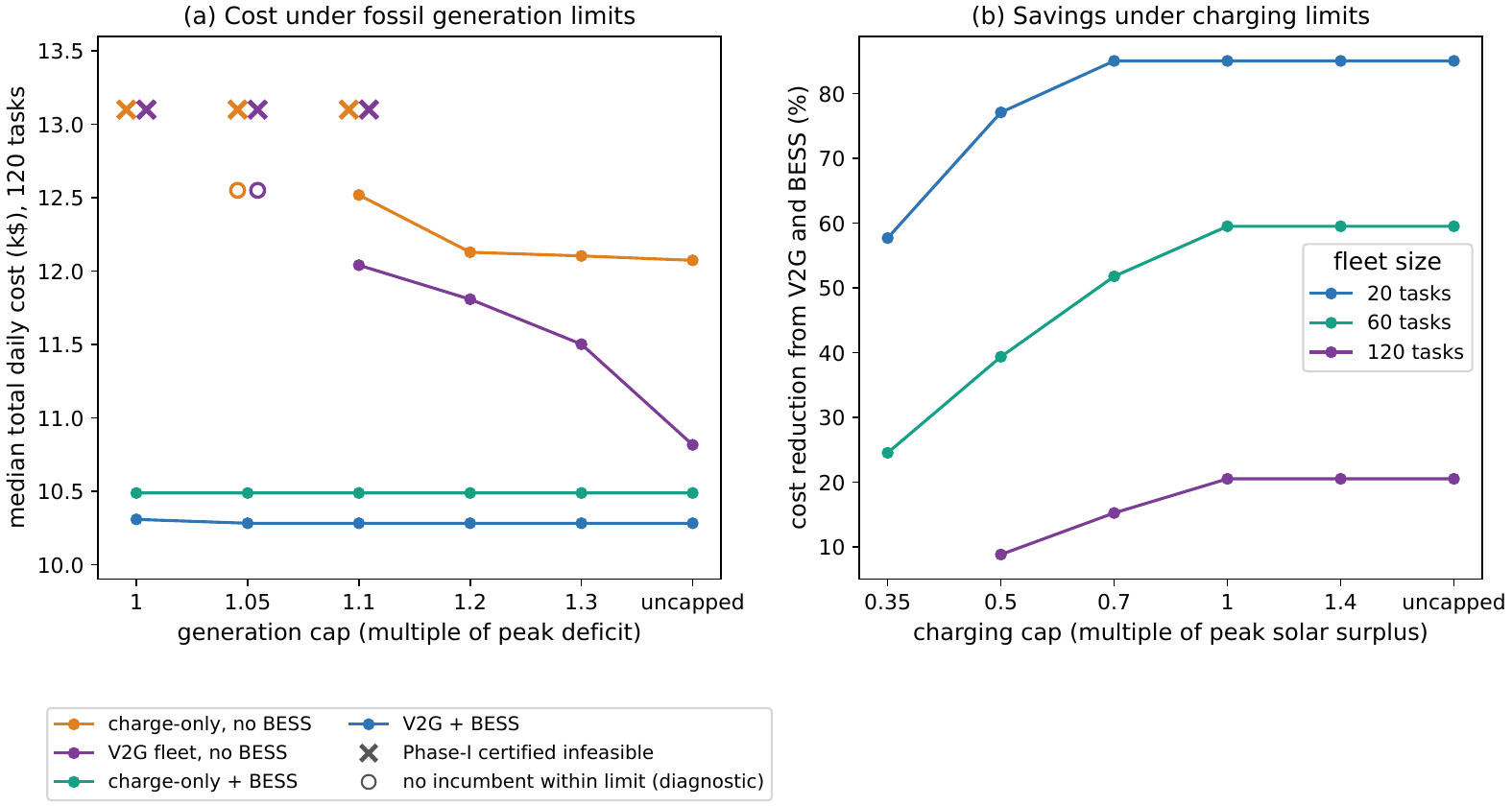}
\caption{Effects of infrastructure limits at $2.5\times$ solar.
Panel (a) reports median cost at 120 tasks as the generation cap changes; the charging cap is fixed at $0.7$ times peak solar surplus.
A cross denotes LP infeasibility certified by Phase I, and an open circle denotes that no feasible integer incumbent was found.
Panel (b) reports the cost reduction from V2G with BESS relative to charge-only without BESS as the charging cap changes.
Each point in panel (b) uses instances with feasible incumbents for both configurations.}
\label{fig:capsfrontier}
\end{figure}

Figure~\ref{fig:capsfrontier}(a) shows the generation cap result at 120 tasks.
In these experiments, every deployed truck must serve at least one transport task.
At a generation cap multiplier of $1.0$, both configurations without BESS are certified infeasible for each of the three generated task sets.
Both configurations with BESS remain feasible over the full range of caps.
Thus stationary storage restores feasibility in these instances.
The experiment does not allow parked V2G trucks that serve only as storage.
The detailed outcomes at 20, 60, and 120 tasks are reported in Appendix~\ref{app:capdetail}.

Figure~\ref{fig:capsfrontier}(b) varies the charging cap while leaving generation uncapped.
A tighter charging cap reduces the cost advantage of adding V2G and BESS because every battery, including stationary storage, must charge through the same capacity constraint.
Comparisons are made only for instances in which both technology configurations have feasible incumbents.
Appendix~\ref{app:capdetail} reports the feasibility counts and cap utilization.

The charging location experiment uses maps with 4--15 task locations and compares charging only at the depot with chargers at selected task locations.
Additional charging locations lower the LP bound for both truck technologies.
On the larger maps at $2\times$ solar, the first two expansion steps, through six charging locations, provide about three quarters of the observed reduction.
At $2\times$ solar, the reported integer incumbents also select fewer stationary batteries when more charging locations are available.
Figure~\ref{fig:stations} reports these results, and Appendix~\ref{app:stationdetail} gives the full design.
The cost comparisons use LP bounds because the final integer solves have substantial reported gaps at the time limit; the battery counts come from integer incumbents and are therefore descriptive.
The experiment does not include charger installation costs or a bidirectional charger premium, so it measures the gross value of charging access rather than a net return on charger investment.

\subsection{Robustness: weather and technology}\label{ssec:robust}

\begin{figure}[htbp!]
\centering
\includegraphics[width=0.78\linewidth]{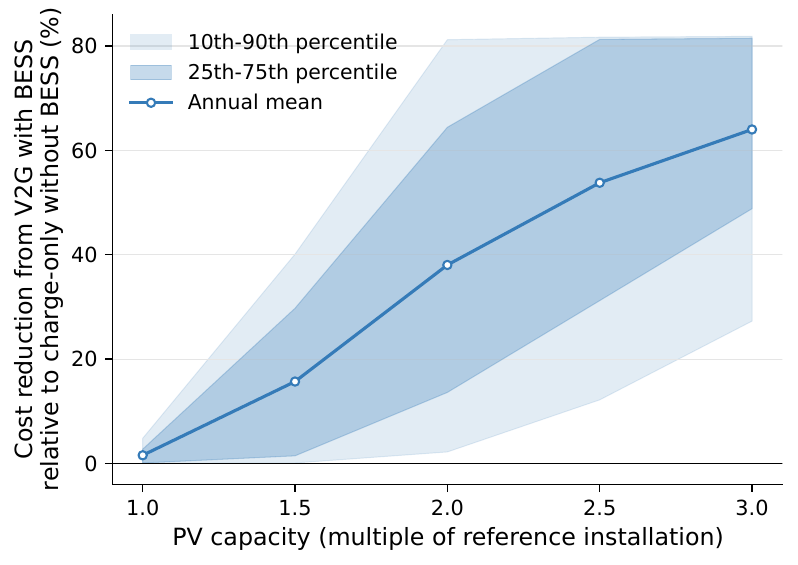}
\caption{Cost reduction from V2G trucks with BESS relative to charge-only trucks without BESS over 365 daily solves using 2023 ERA5 irradiance.
The line is the annual mean; the dark and light bands are the interquartile and interdecile ranges.}
\label{fig:weather}
\end{figure}

\paragraph{A year of weather reanalysis.}
The fixed 60-task benchmark instance is solved for each day of 2023 using ERA5 irradiance at the study location.
Figure~\ref{fig:weather} shows the distribution of the cost reduction from adding V2G and BESS to charge-only trucks without BESS.
The reduction is small at the baseline solar capacity and increases as photovoltaic capacity grows.

Table~\ref{tab:weatherdetail} in Appendix~\ref{app:weatherdetail} reports the median combined reduction and the cost reduction from enabling V2G when BESS is already available.
Table~\ref{tab:climatedetail} reports the corresponding results for four other locations.

At $2.5\times$ solar, the V2G with BESS configuration selects fewer stationary batteries than charge-only with BESS on 345 of 365 days.
It uses slightly more fossil generation on approximately half of the days.
Thus part of its cost reduction comes from replacing stationary batteries with truck storage rather than from additional fossil displacement.
Appendix~\ref{app:weatherdetail} reports these daily comparisons.

Because V2G carries no additional fixed cost in this experiment, allowing discharge cannot raise the optimal cost relative to charge-only operation with BESS.
Weather affects the size of the cost reduction, but not its sign.

\paragraph{Pack size.}
Without BESS, a larger truck battery gives the fleet more capacity to move energy between periods.
Figure~\ref{fig:pack} in Appendix~\ref{app:extra} shows that the estimated cost reduction from V2G rises as pack capacity increases and begins to level off near 1\,MWh.
BESS is unavailable in this experiment, which focuses on truck battery capacity.

\paragraph{Electrification itself.}
For a fixed task set, define the efficiency ratio as the traction energy used by the internal combustion engine (ICE) fleet divided by that used by the electric fleet.
The ICE fuel requirement is linear in this ratio, so equating the two fleet costs gives the break-even ratio analytically.
Across 240 sampled fleets, the largest break-even ratio is 1.97.
The experiment uses ratios from 2.5 to 3.5 as a sensitivity range, not as a universal empirical range.
For context, EPA MOVES4 uses ICE-to-electric energy ratios of 2.6 for combination short-haul trucks and 3.5 for single-unit short-haul trucks, but 2.0 for long-haul trucks; these ratios exclude EV charging losses \citep{epa2023moves4}.
The computed threshold therefore lets a planner insert an efficiency ratio specific to the site and duty cycle rather than assume that electrification always pays.
This comparison adjusts energy use for differences in efficiency; assuming equal energy use for the two technologies reverses the result.
Appendix~\ref{app:electrifydetail} reports summary statistics for the 240 sampled fleets.

\subsection{Resilience: contingencies and endurance}\label{ssec:resil}

The generation capacity experiment in Section~\ref{ssec:infra} chooses assets separately at each capacity level.
Here, a solution for normal operating conditions fixes the number of stationary batteries and sets an upper bound on the trucks available during an outage.
The model may use fewer trucks and reschedule operations.

\begin{figure}[htbp!]
\centering
\includegraphics[width=\linewidth]{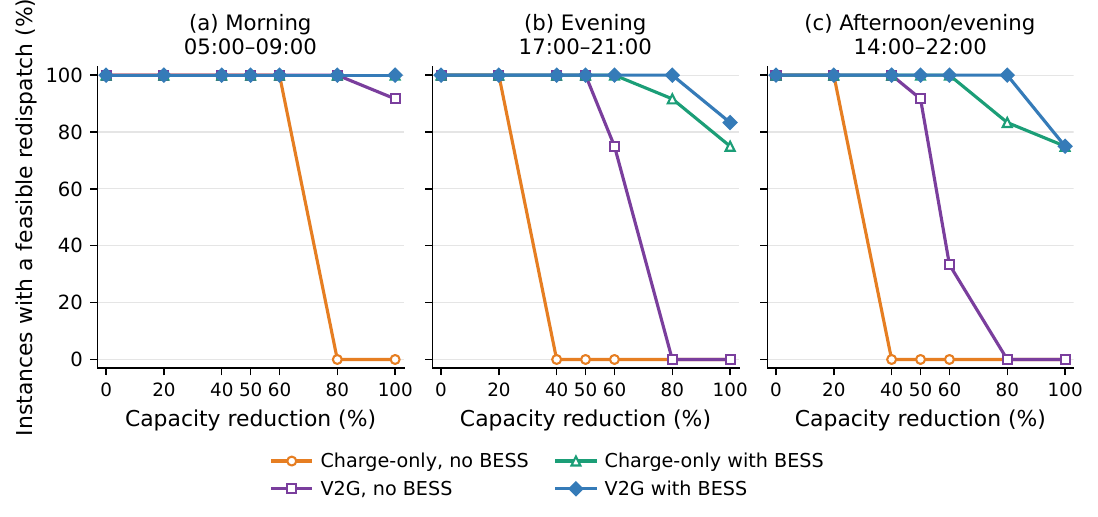}
\caption{Fraction of the 12 instances in each outage window for which a feasible redispatch is found using the asset limits chosen for normal operating conditions.
Panel (a) uses the morning window from 05:00 to 09:00, panel (b) the evening window from 17:00 to 21:00, and panel (c) the longer afternoon-and-evening window from 14:00 to 22:00.
Unresolved instances are not classified as infeasible.}
\label{fig:outage2}
\end{figure}

The generator is derated during three windows: a four-hour morning window from 05:00 to 09:00, a four-hour evening window from 17:00 to 21:00, and an eight-hour afternoon-and-evening window from 14:00 to 22:00.
The eight-hour window contains the four-hour evening window and tests the effect of a longer outage.
Figure~\ref{fig:outage2} reports the three windows separately.

Pooling the two evening-containing windows gives 24 cases at each derating level; these results separate the configurations primarily by BESS.
No feasible incumbent is found for charge-only trucks without BESS beyond a 20\% derating.
V2G without BESS retains a feasible incumbent in every instance through a 40\% derating and in 23 of 24 instances at 50\%.
Both configurations with BESS retain feasible incumbents in every instance through a 60\% derating.

The morning window is less restrictive in these instances because the trucks begin the horizon fully charged and net demand is lower in the morning than during the evening peak.
This is a consequence of the boundary condition and demand profile used in the experiment, rather than a general claim that morning outages are less severe.
The figure reports the fraction of instances for which a feasible incumbent was found; unresolved instances are not classified as infeasible.

A second experiment imposes a daily fuel budget.
Figure~\ref{fig:endur} in Appendix~\ref{app:resildetail} reports the smallest tested budget with a feasible incumbent.
An incumbent for V2G with BESS is found below the fuel baseline without trucks when solar surplus is sufficient; no charge-only incumbent is found at those tested budgets.
The complete counts and fuel budget thresholds are also reported there.

Unresolved fuel budget cases are not infeasibility certificates.
In addition, the deterministic model evaluates scheduling feasibility under a specified outage; it does not estimate generator failure probabilities.

\subsection{Committing schedules using a planning solar profile}\label{ssec:commit}

This experiment commits a schedule before solar generation is known.
A \emph{planning solar profile} is used to choose the truck routes, truck charging plans, and number of stationary batteries.
After the realized solar profile is observed, only BESS dispatch and fossil generation may change.
The experiment uses the fixed 60-task instance.

\begin{figure}[htbp!]
\centering
\includegraphics[width=\linewidth]{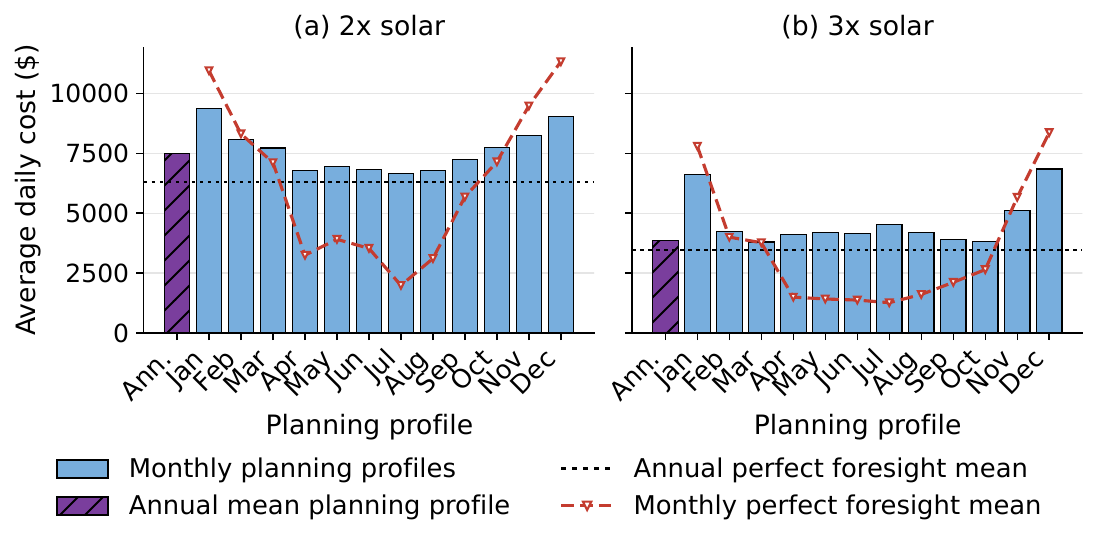}
\caption{Average daily cost in 2023 after committing routes, truck charging, and BESS capacity using one candidate planning solar profile.
BESS dispatch and fossil generation are reoptimized each day.
The dotted line is the annual average of daily perfect foresight incumbents; the dashed line gives the corresponding monthly averages.
Both curves use restricted master incumbents.}
\label{fig:stoch}
\end{figure}

We evaluate thirteen candidate planning profiles on all 365 days of 2023: one based on the annual mean and twelve based on monthly means.
Figure~\ref{fig:stoch} compares their average costs with a daily perfect foresight incumbent benchmark that reoptimizes all decisions each day.
The lowest average cost among the candidate profiles is 5.6\% above that benchmark at $2\times$ solar capacity and 9.7\% above it at $3\times$.

The choice of planning profile matters.
The annual mean profile is 19\% and 12\% above the benchmark at the two solar capacities.
Some winter profiles are 48\% and 98\% above it.
At each solar capacity, the monthly profile with the lowest average cost selects more storage than the annual mean profile.
This storage can be redispatched when the realized solar profile changes.

The planning profiles selected using 2023 are also evaluated on 2022 without reselection.
Every committed schedule remains feasible.
At $2\times$ solar capacity, the 2022 average cost is 5.8\% above the incumbent obtained with perfect foresight and 6.5\% above its LP lower bound.
The corresponding gaps at $3\times$ are 11.4\% and 12.7\%.
Table~\ref{tab:holdout} in Appendix~\ref{app:commitdetail} gives the selected profiles and both comparisons.
This yearlong holdout supports, but does not establish generally, the use of a planning solar profile with operational recourse.

\newpage
\section{Conclusion}\label{sec:concl}
We studied the joint scheduling of electric trucks, V2G operation, and stationary storage in an isolated microgrid.
The model couples truck route selection to the microgrid energy balance.
For a fixed route, charging and discharging are continuous network flow decisions; column generation and dynamic programming handle the exponential route set.

Column generation solved the LP relaxation for every reported scalability instance.
The final schedules come from the integer restricted master and are evaluated against the column generation LP lower bound.

In the comparison of final restricted master incumbents, stationary batteries provide most of the savings from adding storage and V2G: in the median matched instance they account for 96\% of the reduction obtained from the two technologies together.
The additional value of V2G changes with fleet size, solar surplus, task timing, the truck boundary condition, and charger cost.
The estimated break-even value of $\gamma$ changes with charger cost and task duration.
It is therefore best used as a planning indicator that is calibrated to local tasks and costs, rather than as a universal threshold.
Under the tightest tested generation limits, stationary storage also determines feasibility for the largest task sets.

Taken together, the experiments show that V2G should be assessed together with stationary storage, truck availability, and the site's solar profile.
The model can help planners compare investments in trucks, bidirectional chargers, and stationary batteries using local tasks, costs, and energy conditions.
Future work should extend this planning framework to uncertain task demand and solar generation, with schedules that can adapt as conditions are observed.


\begin{appendix}
\section{Proofs and a residual dispatch result}\label{app:proofs}

\paragraph{Proof of Proposition~\ref{thm:np}.}
We reduce from 3-Partition.
Given integers $a_1,\dots,a_{3m}$ with $B/4<a_i<B/2$ and $\sum_i a_i=mB$, the problem asks whether they can be partitioned into $m$ triples that each sum to $B$.
It is strongly NP-complete, meaning that it remains hard when all $a_i$ are polynomially bounded.
Build an instance with $T=3m+1$ blocks and one task location at distance zero from the depot.
Task $j$ occupies block $j$ alone at that location and consumes energy $a_j$; the tasks are temporally disjoint and at the same location, so without energy limits one vehicle could serve them all.
Set base demand and solar to zero ($\Delta_t=0$), battery capacity $G=B$, charging capacity $\kappa_t=0$ for $t\le 3m$ and $\kappa_{3m+1}=\infty$ with rate $\rho=B$, stationary storage excluded, discharging disabled, $\sigma=0$, generation uncapped ($\bar g_t=\infty$ for all $t$), and any $\cv>0$, $\cg>0$.
Set $\eta=0$.
A truck starts full, cannot charge before the final block, and must respect $s_t\ge 0$, so the tasks it covers satisfy $\sum a_j\le B$; conversely any subset with $\sum a_j\le B$ is coverable, with the full recharge boundary restored in the final block.

The full recharge boundary makes every feasible solution recharge exactly $\sum_i a_i=mB$ in total.
With optimal generation dispatch, the fossil term is therefore $\cg\,mB$, and the optimum is $\cg\,mB+\cv\cdot(\text{minimum number of trucks})$.

The minimum truck count is the minimum number of bins of capacity $B$ for the items $a_i$.
It equals $m$ if and only if the 3-Partition instance is feasible: each part then sums to exactly $B$, and the range condition forces three items per part.
Deciding whether the optimum is at most $\cg\,mB+\cv m$ therefore decides 3-Partition; the reduction is polynomial with polynomially bounded data, and the decision version is in NP, so the problem is strongly NP-hard.\qed

For Proposition~\ref{prop:flow}, use the time-expanded network described in Section~\ref{sec:structure}.
Carryover arcs impose battery capacity, charging and discharging arcs are available only while the asset is connected, and fixed withdrawals represent task and deadhead energy.
Charging arcs have gain $1-\eta$; all other energy arcs have unit gain.
The system-level quantities $\Delta_t^-$ and $\Delta_t^+$ enter the master problem, not the feasibility network for one column.

\paragraph{Proof of Proposition~\ref{prop:flow}.}
For a fixed candidate column skeleton, the nonnegative flows on the charge, discharge, and carryover arcs determine the signed energy profile indexed by station and time, $(e_{rht})_{h\in\Hcal^0,\ t\in\Tcal}$.
Energy entering storage from the grid is positive, energy returned from storage to the grid is negative, and no grid interaction gives zero.
Flow conservation with gains along the storage chain is exactly the SoC recursion, including charging loss through the factor $1-\eta$ and fixed withdrawals for task and deadhead traction energy.
The rate limits, battery capacity limits, connection restrictions, and initial/terminal SoC requirements are represented by arc bounds, node bounds, or fixed balances.
Therefore the candidate skeleton satisfies its energy constraints exactly when this time-expanded energy flow LP is feasible.
Applying the cancellation in Lemma~\ref{lem:nosim} with zero prices preserves the SoC path and gives a feasible schedule without simultaneous charging and discharging.

When the LP is feasible, any linear price on the signed grid interactions gives a linear cost on the corresponding grid interface arcs.
For the nonnegative prices in~\eqref{eq:rc}, Lemma~\ref{lem:nosim} gives an optimal schedule without simultaneous charging and discharging.
Replacing $\nu_t e_t^+$ by $\nu_t c_t$ preserves the optimum: the terms agree on such a schedule, and $c_t\ge e_t^+$ for every schedule.
The throughput term is linear in charging and discharging flows, so the pricing energy optimum is obtained by solving the generalized minimum cost flow LP.
When $\eta=0$, the network is an ordinary flow network and the constraint matrix is the node-arc incidence matrix, hence totally unimodular; for integral data, every basic feasible solution is integral.
When $\eta\in(0,1)$, the charging arcs have gain $1-\eta$, so the problem is a generalized flow LP in continuous variables.
In both cases, charging and discharging are represented by continuous flows \citep{ahuja1993network}.\qed

The following lemma concerns residual energy redispatch after the routes and a jointly feasible schedule $p^0$ have been fixed.
We assume $\eta=0$ and define the corresponding residual network $N(S)$ below.

Fix a finite set $S$ of candidate column skeletons.
Their parking patterns, task energy requirements, connection locations, rate limits, and boundary states are fixed.
Also assume a jointly feasible reference schedule $p^0$.
Individual feasibility is not enough because the skeletons share solar surplus and station capacity.

Construct $N(S)$ from the capacity left after $p^0$:
\begin{itemize}[leftmargin=1.5em,itemsep=2pt]
\item Each skeleton has a residual storage chain.
Reverse arcs within the storage and transfer network allow $p^0$ to be redispatched.
\item A source supplies the solar surplus $\Delta_t^-$ that remains unused in each block, and the shared station arcs retain their residual capacities.
\item Each deficit block has one forward demand arc.
Its capacity $u^0_t$ is the part of $\Delta_t^+$ not already displaced by $p^0$.
These arcs have no reverse counterparts, so additional displacement is nonnegative.
\end{itemize}
The reference schedule already satisfies traction and the boundary states.
Hence zero additional flow is feasible in $N(S)$, and feasibility is downward closed.
Let $y_t$ be the additional energy that reduces fossil generation in deficit block $t$; it excludes gross discharge and V2V circulation.

\begin{Lem}[Polymatroid structure of fixed residual dispatch]\label{rem:poly}
When $\eta=0$, let $r_S(B)$ be the maximum additional displacement energy deliverable into blocks $B\subseteq\Tcal$ in the fixed residual network $N(S)$, relative to $p^0$.
Then $r_S$ is normalized, monotone, and submodular, and the achievable additional useful displacement vectors are exactly
\[
\Big\{y\ge 0:\ \sum_{t\in B}y_t\le r_S(B)\ \forall B\subseteq\Tcal,\ 0\le y_t\le u^0_t\Big\},
\]
a polymatroid intersected with the box $[0,u^0]$.
Because the demand arc capacities are part of $N(S)$, the singleton inequalities already imply the box bounds; we display them for clarity.
\end{Lem}

\paragraph{Proof of Lemma~\ref{rem:poly}.}
Collect the demand arcs already present in $N(S)$ at a common super-sink $\psi$.
Let $d_t$ be the tail node of the demand arc for block $t$; this arc, with capacity $u^0_t$, now runs from $d_t$ to $\psi$.
No new arc or capacity is introduced.
For $B\subseteq\Tcal$, let $N_B$ retain only the demand arcs for $t\in B$.
Then $r_S(B)$ is the maximum flow from the source to $\psi$ in $N_B$.
Normalization and monotonicity are immediate: $r_S(\emptyset)=0$, and enlarging $B$ only adds demand arcs.
For submodularity, write every node set $A$ on the source side of a cut (sources in $A$, $\psi\notin A$) as a cut whose value in $N_B$ is $f_A(B)=c_0(A)+\sum_{t\in B}u^0_t\,\mathbf 1[d_t\in A]$, where $c_0(A)$ is the capacity of the other arcs of $N(S)$ that leave $A$.
The demand arcs are accounted for by the displayed sum and are not counted twice.
By the max flow/min cut theorem, $r_S(B)=\min_A f_A(B)$.
Take optimal cuts $A$ for $B$ and $A'$ for $B'$.
The directed cut function $c_0$ is submodular, so $c_0(A\cup A')+c_0(A\cap A')\le c_0(A)+c_0(A')$.
For every $t$, the indicator terms satisfy $\mathbf 1[d_t\in A\cap A']+\mathbf 1[d_t\in A\cup A'] =\mathbf 1[d_t\in A]+\mathbf 1[d_t\in A']$.
This gives equality for $t\in B\cap B'$.
For $t\in B\setminus B'$ (respectively, $t\in B'\setminus B$), the indicator of $A\cap A'$ is bounded by the indicator of $A$ (respectively, $A'$).
Hence
\[
\begin{aligned}
r_S(B\cup B')+r_S(B\cap B')
&\le f_{A\cap A'}(B\cup B')+f_{A\cup A'}(B\cap B')\\
&\le f_A(B)+f_{A'}(B')\\
&=r_S(B)+r_S(B').
\end{aligned}
\]
Finally, a displacement vector $y$ is achievable exactly when the flow with fixed sink demands $y_t$ is feasible.
By the standard supply-demand feasibility theorem, equivalently by adding a supersink and applying max flow/min cut \citep{ahuja1993network}, this holds if and only if $\sum_{t\in B}y_t\le r_S(B)$ for every $B\subseteq\Tcal$, together with the demand arc bounds $y_t\le u^0_t$.\qed

\section{Additional computational results}\label{app:extra}

\subsection{Experiment design and parameter values}\label{app:designdetail}

The reference instance has 60 two-hour tasks over three locations, doubled solar capacity, 700\,kWh truck batteries, a 350\,kW charging rate, no charging loss, and full truck recharge at the end of the horizon.
Fossil generation costs \$0.40/kWh, each truck costs \$45/day, and each 700\,kWh stationary battery costs \$36/day.
Table~\ref{tab:design} lists the departures from this reference instance.

For the studies that share routes, as identified in Section~\ref{ssec:setup}, we run column generation separately for all four technology configurations and combine the generated truck routes before the final integer solves.
The two charge-only cases share the routes generated in those two runs.
The two V2G cases share the routes generated in all four runs, so they also include every charge-only route.
Solving each LP again on its shared route set recovers the column generation bound to numerical tolerance.
This construction removes differences in route generation within each truck technology and ensures that the V2G sets contain every charge-only route.
It does not prove integer optimality over all routes in $\mathcal C$.

Some sensitivity studies, including the battery capacity experiment, use generated route sets constructed separately for each scenario.
Their incumbents are feasible upper bounds, but a difference between incumbents may partly reflect the available routes.
We describe these differences as computational estimates unless the LP and incumbent bounds below certify their sign.

For a matched pair, let $D=z_{\mathrm{base}}-z_{\mathrm{flex}}$ be the absolute daily cost reduction.
If $\mathrm{LB}$ and $\mathrm{UB}$ are the LP lower bound and incumbent upper bound, respectively, then $D\in[\mathrm{LB}_{\mathrm{base}}-\mathrm{UB}_{\mathrm{flex}},\, \mathrm{UB}_{\mathrm{base}}-\mathrm{LB}_{\mathrm{flex}}]$.
The sign of $D$ is certified when this interval excludes zero.

\begin{table}[htbp]
\centering
\caption{Experiment design.
The principal experiments use a 25\,kWh SoC step; experiments using 50\,kWh are identified below.}
\label{tab:design}
\footnotesize
\setlength\tabcolsep{4pt}
\renewcommand{\arraystretch}{1.08}
\begin{tabular}{p{0.21\linewidth}p{0.24\linewidth}p{0.30\linewidth}p{0.13\linewidth}}
\toprule
factor varied & values & departure from reference & result \\
\midrule
workload $n$ & 4--200 tasks & one curve per solar level & Fig.~\ref{fig:boundary} \\
solar capacity & 1--4$\times$ & reference workloads & Fig.~\ref{fig:regime}, Table~\ref{tab:export} \\
truck boundary & free start diagnostic, full recharge, periodic & separate generated columns & \S\ref{ssec:repro}; Table~\ref{tab:boundaryperiodic} \\
timetable & five families & 40/80 tasks, $1.5\times$/$2.5\times$ solar & \S\ref{app:gammadetail} \\
SoC step $\delta$ & 50/25/12.5\,kWh; 6.25\,kWh for the diagnostic with positive charging loss & 60--200 tasks, $2.5\times$ solar, LP & text below \\
pack size $G$ & 350--1{,}400\,kWh & 8--200 tasks & Fig.~\ref{fig:pack} \\
charger power & 100--500\,kW & $G=700$\,kWh & \S\ref{app:gammadetail} \\
generation cap & 1.0--$\infty$ times peak deficit & 20--120 tasks, $2.5\times$ solar & Fig.~\ref{fig:capsfrontier} \\
technology configuration & four configurations, 27 parameter settings, 10 task draws & common truck column construction & Table~\ref{tab:decomp} \\
V2G premium given BESS & \$0/4/8 per truck per day, $\gamma=0.5$--2 & matched solar capacity & \S\ref{app:gammadetail} \\
combined technology premium & \$0/4/8 per truck per day, $\gamma=0.1$--2 & matched solar capacity & \S\ref{app:gammadetail} \\
charging cap & 0.35--$\infty$ times peak solar surplus & generation uncapped & Fig.~\ref{fig:capsfrontier} \\
charging locations & depot, two or five selected task locations, or every task location, $L=4$--15 & a partial set is tested on the larger maps & Fig.~\ref{fig:stations} \\
weather day & 365 days of 2023 & fixed 60-task instance & Fig.~\ref{fig:weather} \\
location & five weather locations & same array and duty & Table~\ref{tab:climatedetail} \\
outage derating & 0, 20, 40, 50, 60, 80, 100\%; three time windows & asset limits for normal operations & Fig.~\ref{fig:outage2} \\
fuel budget & 0.05--2.5 times baseline & 20/60 tasks, $1.5\times$/$2.5\times$ solar & Fig.~\ref{fig:endur} \\
planning solar profile & 13 candidates & fixed 60-task instance & Fig.~\ref{fig:stoch} \\
\bottomrule
\end{tabular}
\end{table}

Refining the SoC step from 50 to 25\,kWh lowers the LP bound by 0.65--0.8\% on the benchmark instances and by at most 2.2\% on the 56-instance alignment study.
Section~\ref{sec:cg} reports the finer comparisons at 12.5 and 6.25\,kWh.

\begin{table}[htbp!]
\centering
\caption{Parameter values, sensitivity ranges, and empirical context for the planning instances.}
\label{tab:prov}
\small
\setlength\tabcolsep{4pt}
\begin{tabular}{p{0.26\linewidth}p{0.2\linewidth}p{0.47\linewidth}}
\toprule
parameter & range used & source or calibration \\
\midrule
fossil generation cost & \$0.20--1.00/kWh & sensitivity range; an estimate for a military microgrid supplied only by diesel is \$0.25/kWh \citep{nrel2011fortcarson} \\
electric truck daily cost & 1.0--2.0 times ICE & scenario multiplier; daily asset cost depends on financing, lifetime, residual value, and utilization; \citet{basma2023total} report Class 8 purchase and ownership cost projections \\
ICE-to-electric energy ratio & 2.5--3.5 & short-haul sensitivity range; EPA MOVES4 uses 2.6 for combination short-haul trucks and 3.5 for single-unit short-haul trucks, but 2.0 for long-haul trucks \citep{epa2023moves4} \\
stationary battery cost & \$26--51/day per 700\,kWh & straight-line allocation of \$200--400/kWh over 15 years, excluding financing, O\&M, augmentation, and salvage; \citet{cole2023storage} project \$245--403/kWh for utility-scale four-hour storage in 2030 \\
bidirectional charger premium & \$0--8 per truck per day & sensitivity range; published V2G equipment costs vary by architecture and power level \citep{steward2017critical} \\
storage discharge cost & \$0--0.13/kWh discharged & stress test applied to truck and stationary battery discharge; the upper value is close in magnitude to a EUR~0.132/kWh mean V2G compensation estimate that includes degradation and infrastructure \citep{sagaria2025vehicle} \\
solar irradiance & 365 days in 2023 & ERA5 reanalysis \citep{hersbach2020era5} \\
charger power & 100--500\,kW & scenario range within the reported range of DC fast charging equipment \citep{afdc2024charging} \\
task energy & 50--300\,kWh & synthetic range from short to energy-intensive tasks \\
\bottomrule
\end{tabular}
\end{table}

\subsection{Boundary condition details}\label{app:boundarydetail}

\begin{table}[htbp]
\centering
\caption{Median increase in final restricted master incumbent cost relative to periodic truck operation with an optimized repeated level over 18 matched bases.
Each fixed level $s_0=s_T=c$ applies to every truck; stationary storage remains periodic.}
\label{tab:boundaryperiodic}
\footnotesize
\begin{tabular}{lccccc}
\toprule
 & $c=G$ & $c=3G/4$ & $c=G/2$ & $c=G/4$ & $c=0$ \\
\midrule
charge-only without BESS & 14.7\% & 7.8\% & 3.2\% & 6.5\% & 13.4\% \\
charge-only with BESS    & 2.8\%  & 1.8\% & 1.1\% & 2.0\% & 2.8\%  \\
V2G without BESS         & 28.1\% & 6.7\% & 0.5\% & 1.2\% & 18.8\% \\
V2G with BESS            & 0.1\%  & 0.0\% & 0.0\% & 0.0\% & 0.0\%  \\
\bottomrule
\end{tabular}
\end{table}

All 432 boundary solves produce feasible incumbents.
Every column generation LP terminates with no route whose reduced cost is below the stopping threshold.
Of the final integer solves, 43 reach the time limit.
The LP bound with an optimized repeated level is no greater than the bound for the corresponding fixed level in every comparison.
Fixing $c=G/2$ reduces the median column generation time by a factor of 6.1.
For the configurations without BESS, intermediate repeated levels are least costly: low levels leave less energy for morning tasks, while high levels require more evening recharge.

\begin{table}[htbp!]
\centering
\caption{Change in fossil use relative to the microgrid without trucks, in equivalent gallons per day at 33\,kWh per gallon.
Negative entries denote displacement of base fossil generation.
A star indicates zero total fossil generation.
Each row uses one task set and full truck recharge.}
\label{tab:export}
\small
\setlength\tabcolsep{4.5pt}
\begin{tabular}{rrrrrrrr}
\toprule
tasks $\backslash$ surplus & 3.7 & 10.3 & 16.8 & 23.7 & 30.9 & 38.0 & 45.1 \,MWh/d \\
\midrule
20  & $+48$   & $-148$ & $-348$ & $-450$     & $-446^{*}$ & $-430^{*}$ & $-430^{*}$ \\
40  & $+211$  & $+18$  & $-182$ & $-388$     & $-446^{*}$ & $-430^{*}$ & $-430^{*}$ \\
60  & $+373$  & $+173$ & $-24$  & $-233$     & $-446^{*}$ & $-430^{*}$ & $-430^{*}$ \\
80  & $+508$  & $+311$ & $+120$ & $-94$      & $-308$     & $-430^{*}$ & $-430^{*}$ \\
100 & $+685$  & $+485$ & $+289$ & $+80$      & $-138$     & $-353$     & $-430^{*}$ \\
120 & $+823$  & $+622$ & $+429$ & $+215$     & $+0$       & $-218$     & $-418$ \\
140 & $+1016$ & $+809$ & $+618$ & $+418$     & $+194$     & $-24$      & $-236$ \\
\midrule
max. displacement & $-512$ & $-488$ & $-461$ & $-452$ & $-446$ & $-430$ & $-430$ \\
\bottomrule
\end{tabular}
\end{table}

The maximum displacement row becomes less negative as solar capacity increases because the base itself burns less fossil fuel, leaving less generation for the fleet to displace.

\subsection{Additional results for $\gamma$}\label{app:gammadetail}

The ratio $\gamma$ does not by itself determine the cost reduction when the timing of tasks changes.
For example, consider two 200\,kWh tasks and one two-hour 800\,kWh solar surplus period, with BESS unavailable.
Scheduling both tasks during that period gives 22.6\,MWh of fossil generation.
Scheduling them outside the surplus period gives 22.0\,MWh with charge-only trucks and 21.8\,MWh with V2G trucks.
All four integer restricted masters are solved to optimality over their sets of generated columns.
Across five timetable families, retiming the same tasks materially changes fossil generation at fixed $\gamma$.
In a separate sensitivity experiment with stationary storage available and charging uncapped, setting deadhead energy and truck deployment cost to zero eliminates this timetable spread.
This indicates that the spread is driven by deadhead energy, truck deployment cost, or their interaction, none of which is represented by $\gamma$.

The bidirectional charger premium is varied from \$0 to \$8 per truck per day as a sensitivity analysis.
Published V2G equipment costs vary substantially with power level and whether the bidirectional conversion equipment is on the vehicle or in the charger \citep{steward2017critical}; the tested daily premium is therefore not presented as a universal hardware cost estimate.
In the 60-task experiment that scales solar capacity, we sample $\gamma$ from 0.10 to 0.80 in steps of 0.05 for five task sets.
At the \$8 premium, the median saving changes sign between $\gamma=0.40$ and $0.45$, and the crossings across the task sets range from 0.37 to 0.48.
The saving is certified negative for every task set through 0.35 and positive for every task set from 0.55.
At the \$4 premium the median crossing is near $\gamma=0.35$.
With no premium, the median saving is nonnegative at every tested value.
In a separate sensitivity test with a coarser SoC step, an additional cost of \$0.13/kWh on discharge from both trucks and stationary batteries shifted the estimated break-even upward by about 0.10 at the \$8 premium.

The comparison between V2G and charge-only trucks after BESS is available comprises 630 solves: seven task sets, three workloads, values of $\gamma$ from 0.5 to 2, and premiums of \$0, \$4, and \$8 per truck per day.
At the \$8 premium, the 20- and 60-task incumbent medians are negative at every sampled $\gamma$, and the 120-task medians turn positive near $\gamma=1$, reaching 5.2\% at $\gamma=2$.
At the \$4 premium, the 60-task crossing is near $\gamma=1$ and the 120-task crossing lies below $\gamma=0.75$, reaching 9.4\% at $\gamma=2$.
Of the 105 intervals for the matched cost difference at the \$8 premium, 21 certify a positive effect, 70 certify a negative effect, and 14 contain zero.
For those 14 intervals, the sign of the incumbent difference is descriptive rather than certified.

For the comparison that adds both V2G and BESS to charge-only trucks without BESS, the fixed base experiment uses the same paired route set construction and gives a range over which the median changes sign rather than a single crossing.
At the \$8 premium, the medians for all three fleet sizes are negative through $\gamma=0.35$ (from about $-0.1\%$ at 20 tasks to $-1.3\%$ at 120 tasks) and positive at $\gamma=0.5$ ($+0.2\%$, $+1.1\%$, and $+1.4\%$; six of seven task sets positive at 20 tasks).
Above $\gamma=0.5$ the medians grow with fleet size, reaching $+21\%$, $+58\%$, and $+77\%$ at $\gamma=2$.
At $\gamma\approx0.9$, adding both technologies reduces cost by 4.6\%, 10.0\%, and 23.9\% for 20, 60, and 120 tasks, respectively.
At a \$4 premium, the 20-task median turns positive near $\gamma=0.3$.
These fixed base results explain why the crossing observed when the fleet and microgrid are scaled together should not be treated as universal.

Shifting the solar or demand profile by up to two hours while preserving daily totals changes the combined technology cost reduction by less than 0.8 percentage points in this experiment.
Reducing the charge rate to 100\,kW materially lowers the reduction in the middle of the range; the 200 and 350\,kW results are similar.

\begin{figure}[htbp!]
\centering
\includegraphics[width=\linewidth]{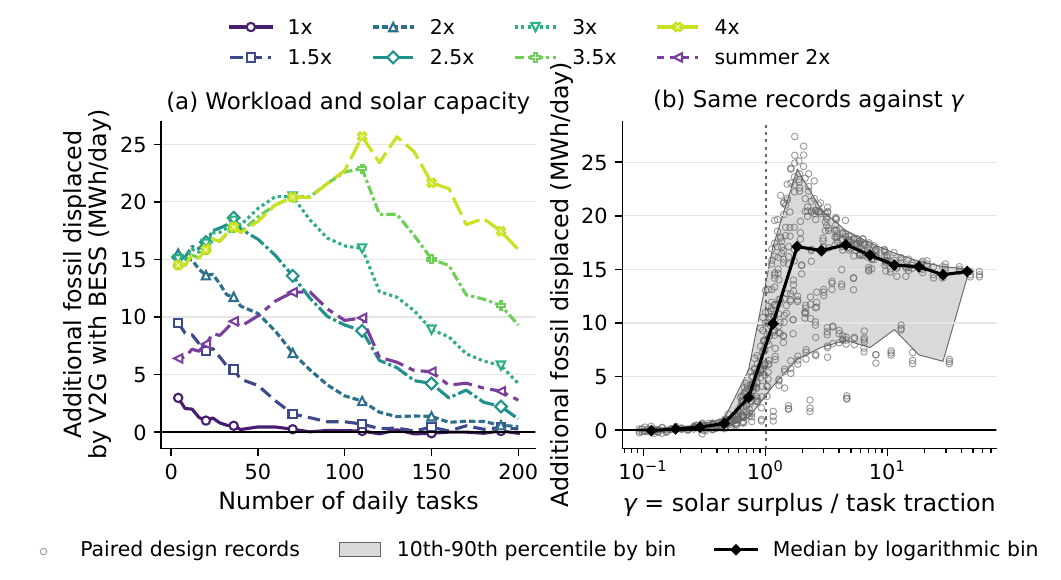}
\caption{Additional fossil generation displaced by V2G trucks with BESS relative to charge-only trucks without BESS.
The left panel varies workload at eight solar profiles; points are medians over three task draws.
The right panel replots the same 624 paired records against $\gamma$ and shows the median and 10th--90th percentile band within each bin.
This figure concerns fossil displacement, not the cost break-even values reported in Section~\ref{ssec:when}.
Small negative differences are descriptive: the plotted schedules minimize total cost rather than fossil use and are not proven integer optima.}
\label{fig:boundary}
\end{figure}

\subsection{Infrastructure details}\label{app:capdetail}

For the 120-task instances certified infeasible, the smallest positive Phase I objective is about $0.047$.
Here $\epsilon=10^{-6}$ and $n\epsilon=1.2\times10^{-4}$.
The recorded value therefore exceeds the sufficient certification threshold by more than two orders of magnitude.
These runs use $\eta=0$, full truck recharge, and the same family of route skeletons and divisibility conditions used in the main capacity experiment.

At 120 tasks and a generation cap equal to the base's peak net demand before trucks or stationary batteries are added, all three task sets are certified infeasible by Phase I for both configurations without BESS.
At 1.05 times the peak, two task sets are certified infeasible and one is unresolved.
At 1.1 times the peak, one task set is certified infeasible, one is unresolved, and one is feasible.
Most task sets are feasible from 1.2 times the peak.
At 1.3 times the peak, the median costs are \$12{,}103 (charge-only without BESS), \$11{,}501 (V2G without BESS), \$10{,}487 (charge-only with BESS), and \$10{,}281 (V2G with BESS).
All four configurations are feasible at 20 and 60 tasks throughout this cap range.
Across those cases, charge-only and V2G without BESS cost as much as $6.0$ and $3.9$ times the cost of charge-only operation with BESS at 20 tasks; at 60 tasks, the corresponding ratios are $2.0$ and $1.45$.

With generation uncapped, tightening the charging cap from no cap to $0.35$ times peak solar surplus lowers the cost reduction from adding V2G and BESS from 85\% to 58\% at 20 tasks and from 60\% to 25\% at 60 tasks.
At 120 tasks, the corresponding reduction falls from 21\% without a cap to 9\% at $0.5$ times peak solar surplus; the $0.35$ case has no matched feasible pair.
The charging cap is fully used in the V2G solutions through 0.7 times peak solar surplus, and peak utilization is 0.71 at 1.4 times peak surplus.
At 0.35 times peak surplus and 120 tasks, eight of the twelve task set/configuration cases are certified infeasible; the other four have feasible LP solutions but no feasible integer incumbent within the time limit.
No matched comparison is reported for that point.

\subsection{Charging location details}\label{app:stationdetail}

The location study uses maps with 4--15 task locations, four solar levels, and 40--160 one-hour tasks.
It uses three task draws for each $1\times$ solar profile and nine for each $2\times$ profile.
Charging only at the depot is compared with a charger at every task location.
On the larger maps, partial sets add chargers at two or five task locations.
The 3{,}456 solves are compared through their column generation LP bounds because the integer restricted masters have substantial reported gaps at the time limit.

\begin{figure}[htbp!]
\centering
\includegraphics[width=\linewidth]{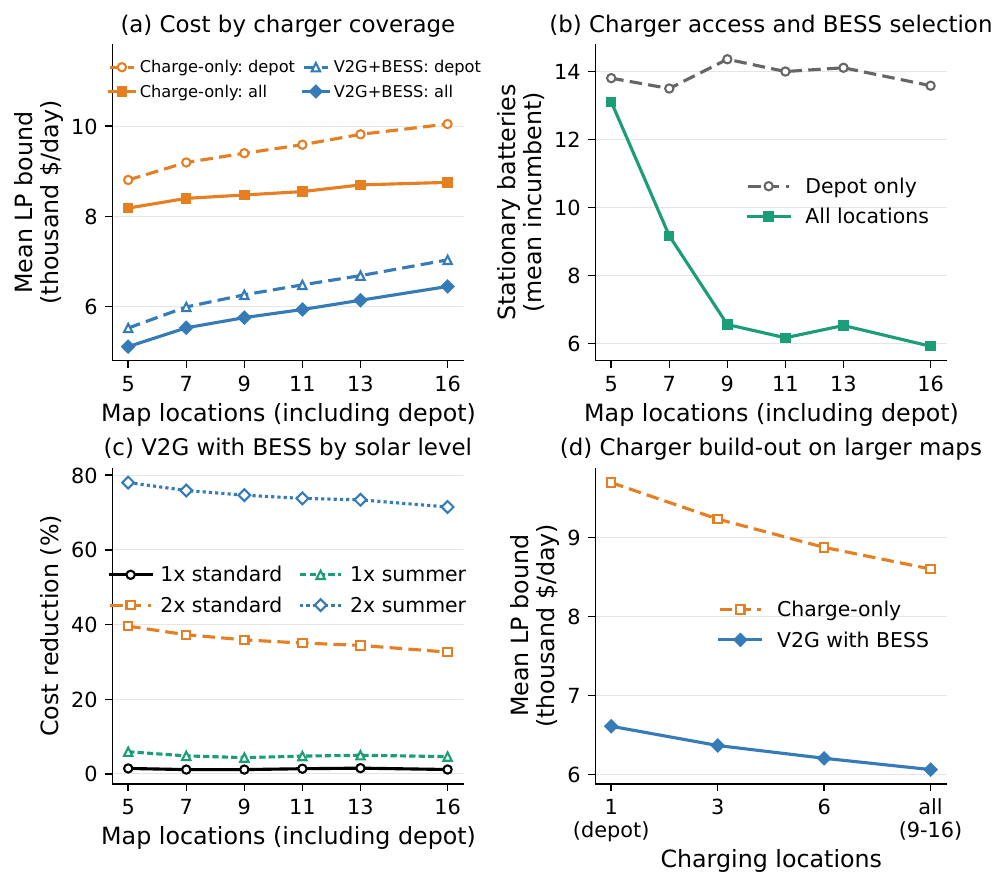}
\caption{Effect of charging locations.
Panel (a) compares charging only at the depot with charging at every task location.
Panel (b) reports stationary batteries selected by V2G with BESS at $2\times$ solar.
Panel (c) reports the cost reduction from V2G with BESS relative to charge-only trucks without BESS at four solar levels.
In panels (a)--(c), the horizontal coordinate counts all map locations, including the depot; depot-only cases have one charging location.
Panel (d) varies the number of charging locations and uses the six task draws available for every choice of charging locations.
Panels (a), (c), and (d) use LP bounds.
Panel (b) reports stationary battery counts from integer incumbents obtained at the time limit.
No mandatory technology premium is imposed.}
\label{fig:stations}
\end{figure}

\subsection{Weather details}\label{app:weatherdetail}

\begin{table}[htbp]
\centering
\caption{Median percentage cost reductions over 365 daily solves at the study location.
The combined comparison is V2G with BESS versus charge-only without BESS.
The conditional comparison is V2G with BESS versus charge-only with BESS.}
\label{tab:weatherdetail}
\small
\begin{tabular}{rrr}
\toprule
solar capacity & combined comparison & conditional V2G comparison \\
\midrule
$1.0\times$ & negligible & 0.2\% \\
$1.5\times$ & 10\% & 0.8\% \\
$2.0\times$ & 33\% & 2.1\% \\
$2.5\times$ & 57\% & 3.4\% \\
$3.0\times$ & 81\% & 8.1\% \\
\bottomrule
\end{tabular}
\end{table}

At $2.5\times$ solar capacity, the V2G configuration selects fewer stationary batteries than the charge-only configuration on 345 of 365 days.
It uses slightly more fossil generation on approximately half of the days, showing that part of the cost reduction comes from replacing stationary capacity rather than from additional fossil displacement.

\begin{table}[htbp]
\centering
\caption{Weather location sensitivity.
Entries are percentage cost reductions.
The last group spans two sites.}
\label{tab:climatedetail}
\small
\begin{tabular}{lrrp{0.34\linewidth}}
\toprule
location group & median at $2\times$ & median at $3\times$ & distributional result \\
\midrule
Southwest Asia & 50\% & 82\% & 10th percentiles 16\% and 54\%; conditional V2G medians 2.8\% and 8.3\% \\
Northeast Asia & 8\% & 44\% & 10th percentile near zero \\
\shortstack[l]{Northern Europe /\\Arctic Circle} & 0\% & 0\% & 90th percentiles 19--22\% and 74--77\%; conditional median zero \\
\bottomrule
\end{tabular}
\end{table}

Both sites in the Northern Europe and Arctic Circle group have zero median reductions at $2\times$ and $3\times$ solar.
The large upper tails in the final row show that the cost reduction is strongly seasonal.

\subsection{Electrification details}\label{app:electrifydetail}

Across 240 sampled fleets, the mean break-even efficiency ratios are 1.33, 1.47, and 1.60 when the electric truck daily cost is 1.0, 1.5, and 2.0 times the ICE truck cost.
The corresponding 95th percentiles are 1.48, 1.67, and 1.88, and the largest observed value is 1.97.
At \$0.40/kWh, electrifying the 60-task reference fleet lowers modeled daily operating cost by about \$8{,}500.
The comparison deploys 18 ICE trucks and 30 electric trucks, so this result for the full fleet is not converted into a payback estimate for each vehicle.

\subsection{Outage and fuel budget details}\label{app:resildetail}

Pooling the four-hour evening and eight-hour afternoon-and-evening windows gives 24 cases at each derating level.
V2G without BESS has a feasible incumbent in every pooled case through a 40\% derating and in 23 of 24 cases at 50\%.
Both configurations with BESS have feasible incumbents in every pooled case through 60\%; at total generation loss, the counts are 18 of 24 for charge-only with BESS and 19 of 24 for V2G with BESS.
During the morning window, all four configurations survive through 60\%; at total loss the counts are 0 of 12 for charge-only without BESS, 11 of 12 for V2G without BESS, and 12 of 12 for each configuration with BESS.

At $2.5\times$ solar, the smallest tested fuel budget with a feasible V2G with BESS incumbent is 5\% of the baseline without trucks.
The analytic ICE requirement for that instance is 38 times as large; across the tested instances, this ratio is at least 2.3.
At $1.5\times$ solar and 60 tasks, the smallest tested budget with a feasible incumbent is 1.3 times the baseline; smaller tested budgets remain unresolved.

These runs with a fixed fleet omit the dual $\tau\le0$ associated with the fleet cap from pricing.
The omitted reduced cost term is $-\tau\ge0$.
With this omission, these runs are used only to report feasible incumbents; we do not claim LP optimality or infer infeasibility from unsuccessful searches.

\begin{figure}[htbp!]
\centering
\includegraphics[width=0.75\linewidth]{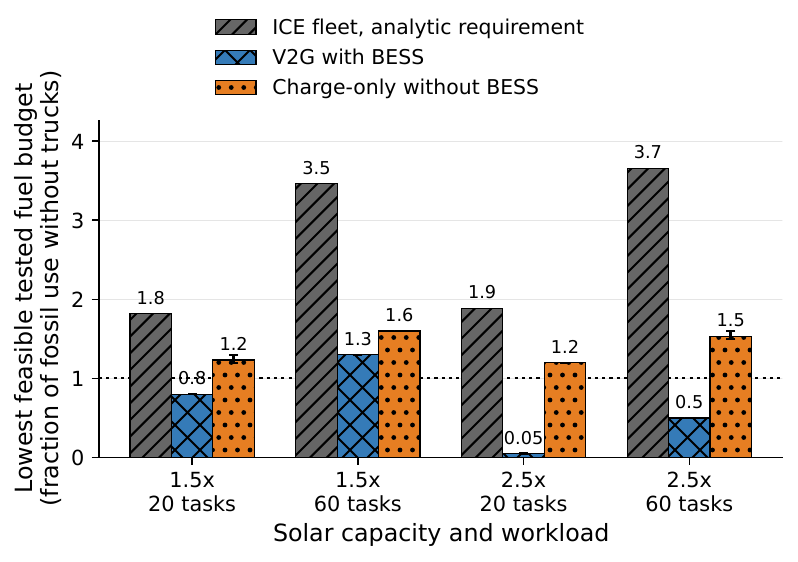}
\caption{Lowest tested daily fuel budget with a feasible incumbent, as a multiple of fossil use in the microgrid without trucks.
Gray bars show the analytic ICE requirement.
Blue and orange bars show the mean over three task sets for V2G with BESS and charge-only trucks without BESS; whiskers span the three values.
Lower tested budgets are unresolved, not certified infeasible.}
\label{fig:endur}
\end{figure}

\subsection{Weather holdout}\label{app:commitdetail}

\begin{table}[htbp]
\centering
\caption{Evaluation on 2022 after selecting the planning solar profile using 2023.
The gap is reported relative first to the daily incumbent obtained with perfect foresight and then to its LP lower bound.}
\label{tab:holdout}
\small
\begin{tabular}{lcccc}
\toprule
solar capacity & 2023 profile & feasible days &
gap to incumbent & gap to LP bound \\
\midrule
$2\times$ & July & 365/365 & 5.8\% & 6.5\% \\
$3\times$ & March & 365/365 & 11.4\% & 12.7\% \\
\bottomrule
\end{tabular}
\end{table}

Winter planning profiles remain the most expensive in the holdout, with gaps of 62--112\%.
The best candidate depends on solar capacity, so the experiment does not support a universal choice of planning month.

\begin{figure}[htbp!]
\centering
\includegraphics[width=\linewidth]{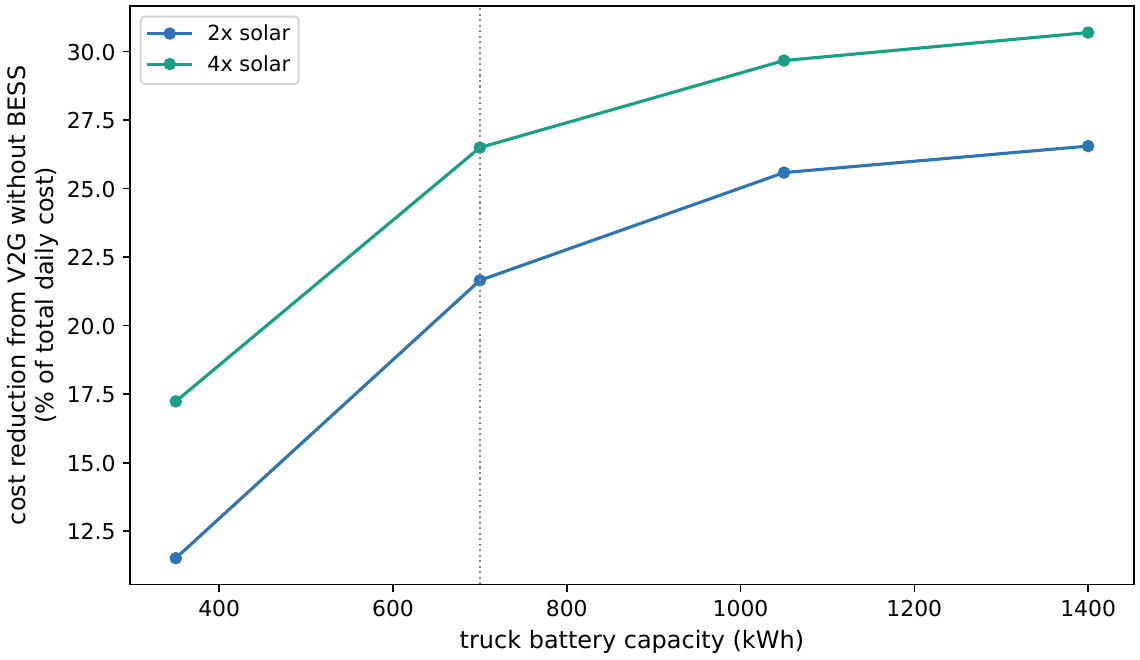}
\caption{Estimated effect of truck battery capacity on the cost reduction from V2G without BESS for fixed task sets.
Each battery capacity scenario uses its own generated route set.
The reduction rises from 12--17\% of daily cost at 350\,kWh to 27--31\% at 1.4\,MWh and begins to level off above approximately 1\,MWh.}
\label{fig:pack}
\end{figure}

\end{appendix}

\section*{Declarations}

\paragraph{Ethics approval and consent to participate}
Not applicable.

\paragraph{Consent for publication}
Not applicable.

\paragraph{Funding}
N.C. and A.S. received funding from the U.S.\ Department of Defense,
Office of Naval Research, under grant N00014-23-1-2349.

\paragraph{Competing interests}
The authors declare no competing interests.

\paragraph{Availability of data and materials}
The research code, input profiles, and selected result files are publicly
available at \url{https://github.com/ndandnd/evspv2g_dp}.

\paragraph{Acknowledgements}
The authors gratefully acknowledge Prof.\ Tong Wu, Department of Electrical
and Computer Engineering, University of Central Florida, for guidance and
for preparing the data used in this study.

\bibliographystyle{spbasic}
\bibliography{references}

\end{document}